\documentclass[a4paper]{amsart}

\usepackage{epsfig,amsmath}
\usepackage{xspace}
\usepackage[psamsfonts]{amssymb}
\usepackage[latin1]{inputenc}
\usepackage{float}\restylefloat{figure}
\usepackage{color}

\newtheorem*{lemma*}{\bf Lemma}
\newtheorem*{sublemma*}{\bf Sublemma}
\newtheorem*{claim*}{\bf Claim}
\newtheorem*{complement*}{\bf Complement}

\renewcommand{\epsilon}{\varepsilon}

\newcommand{\cal}[1]{\mathcal{#1}}
\newcommand{\setof}[2]{\big\{{#1}\,\big|\,{#2}\big\}}

\newcounter{rememberItem}

\def\bEA{\begin{eqnarray*}}
\def\eEA{\end{eqnarray*}}
\def\bEAn{\begin{eqnarray}}
\def\eEAn{\end{eqnarray}}

\def\ds{\displaystyle}
\def\tend{\longrightarrow}

\def\ds{\displaystyle}
\def\wt{\widetilde}

\def\on{\operatorname}

\def\cal{\mathcal}

\def\C{{\mathbb C}}
\def\D{{\mathbb D}}

\def\R{{\mathbb R}}
\def\S{{\mathbb S}}

\def\Z{{\mathbb Z}}

\def\Re{{\on{Re}\,}}
\def\Im{{\on{Im}\,}}

\newtheorem{proposition}{Proposition}
\newtheorem{theorem}[proposition]{Theorem}
\newtheorem{lemma}[proposition]{Lemma}

\newcommand{\sq}{\mathrm{Sq}}
\newcommand{\cotan}{\on{cotan}}

\newcommand{\titre}[1]{\textsl{#1}:\ }

\usepackage{pstricks}

\usepackage[ps2pdf]{hyperref}

\begin{document}


\author{Arnaud Chéritat}
\title[A worked out example of straightening]{Beltrami forms, affine surfaces and the Schwarz-Christoffel formula: a worked out example of straightening}
\begin{abstract}Consider the straightening $\phi$ of a Beltrami form that is constant on a square, with the corresponding ellipses having a vertical major axis, and null outside. A generalized Schwarz-Christoffel formula is used to express the inverse of $\phi$. The formula is found by introducing an affine Riemann surface. This formula is used to draw on a computer the image of the square by $\phi$, and practical aspects are discussed. The resulting shapes are shown for different values of the constant dilatation ratio of the ellipses ($=$major axis$/$minor axis). The limit when this ratio tends to infinity is surprising. A model of this limit is proposed, produced by an affine surface uniformization.
\end{abstract}
\maketitle

In all this article, unless otherwise stated, the word affine means complex-affine. The term ``affine map'' will alway refer to non constant affine maps. The word surface will \emph{not} refer to complex manifolds of complex dimension $2$.

\section{Structure of the article}

Section~\ref{sec:teaser} plays the role of an introduction. There, we present the problem that motivated this work. It involves a Beltrami form that is constant in a square and zero outside. What does the image of the square look like? What is the limit of this shape when the dilatation ratio tends to infinity? The author shows a few early computer experiments.

In section~\ref{sec:sol} we explain how to express the inverse of the straightening of the Beltrami form using a generalized Schwarz-Christoffel formula. The problem of straightening the Beltrami form is transformed into the problem of conformally uniformizing to a punctured Riemann sphere some affine surface. We introduce a differential invariant characterizing the induced affine surface structure on the punctured Riemann sphere; it turns out to be a simple rational map. From this we deduce the formula.

In section~\ref{sec:prac} we discuss practical aspects of numerically using the formula.

In section~\ref{sec:limit} we study the limit of the formula as the dilatation ratio tends to infinity, and give a candidate affine surface whose uniformization would correspond to that limit.

\section{A teaser}\label{sec:teaser}

In this section, we will present two questions, that originally motivated the work of the author, and show a few computer experiments. The tools to answer the second question are introduced in section~\ref{sec:sol}. This will allow us to run better computer experiments in section~\ref{sec:result}, helping us to guess the answer. The same tools allow to give a precise candidate for the answer, and this is discussed in section~\ref{sec:limit}.

In subsection~\ref{subsec:pb} we state the questions.
In subsection~\ref{subsec:exper} we begin exploring it with computer experiments.
In subsection~\ref{subsec:guess} we give a few naïve guesses of the possible answer.

\subsection{Straightening a square}\label{subsec:pb}
Consider the square $\sq$ in $\C$ defined by $x+iy\in \sq$ $\iff$ $x\in[-1,1]$ and $y\in[-1,1]$.
Consider the Beltrami form $\mu(z)\frac{d\bar{z}}{dz}$ with $\mu(z)=0$ outside $\sq$ and $\mu(z)=\frac{1}{3}$ in $\sq$. The corresponding ellipses are circular outside $\sq$ and vertical with ratio $K=2$ in $\sq$. The straightenings of this Beltrami form are quasiconformal over $\C$, and conformal outside the square.
Let $\phi:\C\to\C$ be the unique straightening satisfying the following normalization: 
\[\phi(z) \underset{z\tend\infty}{=} z + 0 +o(1).\]

\begin{figure}[htbp]
\scalebox{0.4}{\includegraphics{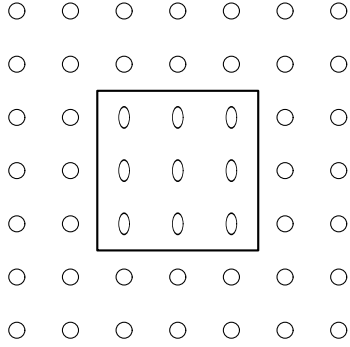}}
\caption{The Beltrami form $\mu$, visualized as a field of ellipses.}
\end{figure}

\medskip

\textbf{First question}: what does the image of the square look like?

\medskip

\begin{figure}[htbp]
\begin{picture}(350,110)(0,0)
\put(-5,0){\scalebox{0.5}{\includegraphics{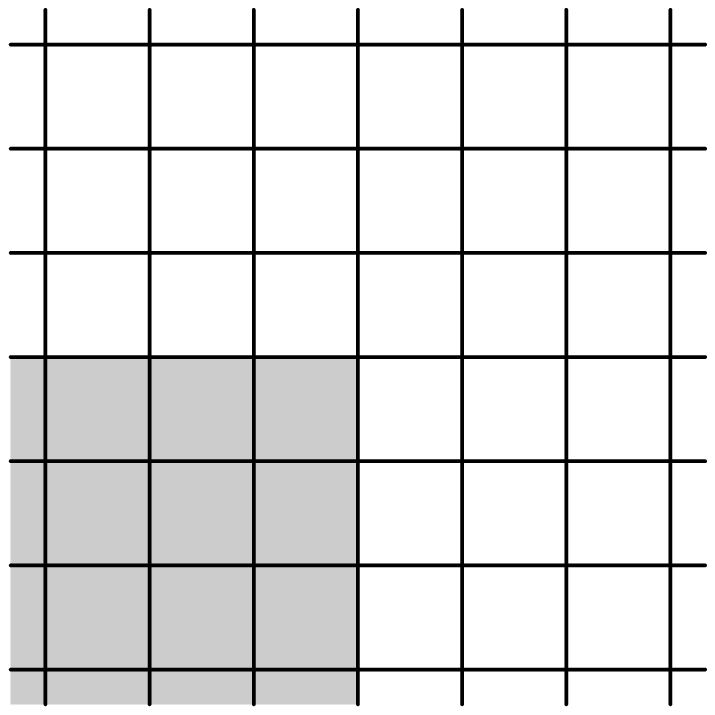}}}
<\put(105,52.5){$\overset{f}{\longrightarrow}$}
\put(120,0){\scalebox{0.5}{\includegraphics{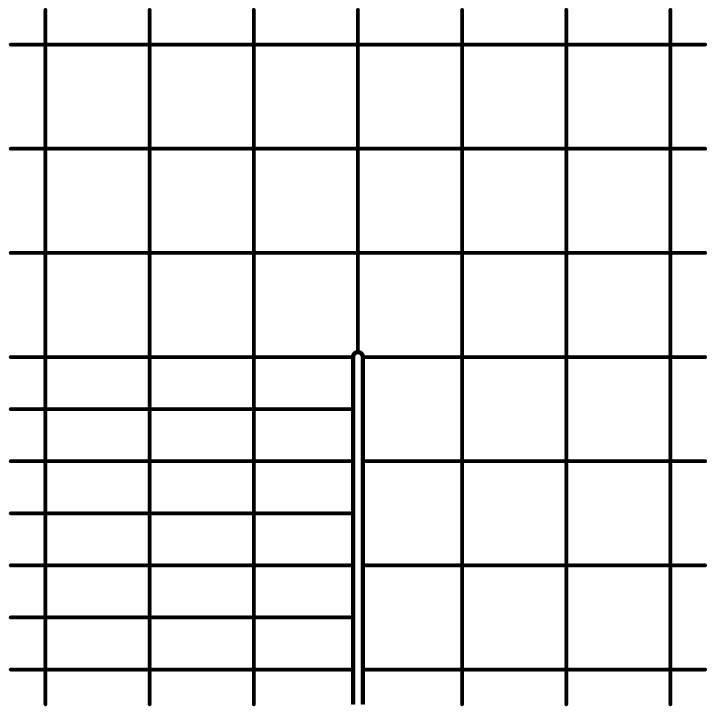}}}
\put(230,52.5){$\overset{z^\alpha}{\longrightarrow}$}
\put(245,0){\scalebox{0.5}{\includegraphics{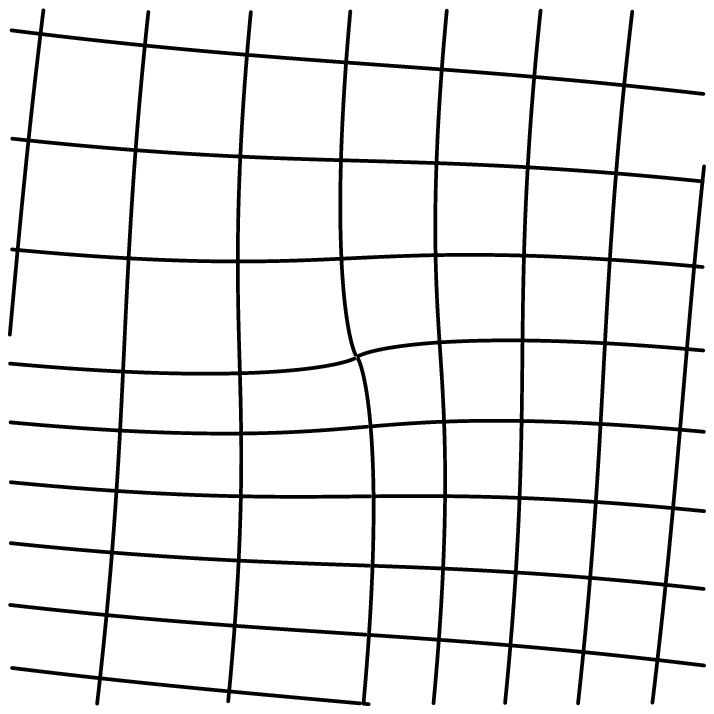}}}
\end{picture}
\caption{Straightening the Beltrami form corresponding to horizontal ellipses with ratio 2 on the gray quadrant and circles outside.}
\label{fig:b}
\end{figure}

\begin{figure}[htbp]
\rotatebox{90}{\includegraphics{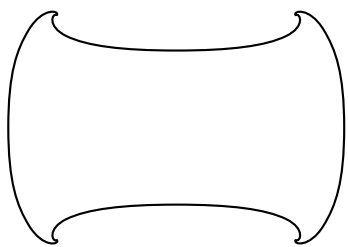}}
\caption{A sketch of what the image of the square is expected to look like, with the spirals exaggerated. It was hand drawn with a Vector Drawing Program.}
\label{fig:below}
\end{figure}

One may give a rough guess, with the following heuristics:
If one zooms on one corner, say the upper right one, $\mu$ tends to some Beltrami form $\nu(z)\frac{d\bar{z}}{dz}$ with $\nu(z)=1/3$ on the lower left quadrant, and $0$ anywhere else. This Beltrami form is explicitly solvable, by the following method, illustrated on figure~\ref{fig:b}: slit the complex plane along the vertical half-line $L=(-\infty,0]i$. Define a map $f:\C\setminus L\to\C\setminus L$ by $f(x+iy) = x+i\frac{y}{2}$ in the lower left quadrant, and $f(z)=z$ elsewhere. This map solves the Beltrami equation, but has a discontinuity along $L$: a point that crosses this boundary from the lower left quadrant to the lower right quadrant sees its image by $f$ jump from $iy/2\in L$ to $iy$.
The idea is then to define a Riemann surface $\cal S$ obtained by gluing the lower left and the lower right quadrants along their common boundary $L\setminus\{0\}$ by identifying $z/2$ to $z$. It is known that this Riemann surface is uniformized to $\C^*$ by a branch of the map $z\mapsto z^\alpha$ for some well chosen complex $\alpha$. Indeed in the coordinates $w=\log z$, the problem is equivalent to gluing the band ``$\Im z \in [-\pi/2,3\pi/2]$'' by $x-i\pi/2\sim x+3i\pi/2-\log 2$, i.e.\ by the translation of vector $2i\pi-\log 2$. The job is done by for instance by $w\mapsto \exp(2 i \pi w/(2i\pi-\log 2))$. So $\alpha=2i\pi/(2i\pi-\log 2)=1/(1-\log(2)/2i\pi)$.
The image of the two half-lines composing the boundary of the lower left quadrant (and of any half-line through $0$) are logarithmic spirals. They turn quite slowly: to make one turn, the distance to $0$ needs to be divided by $e^{(2\pi)^2/\log 2}\approx 6\times10^{24}$.
So we may expect the image of the square to look roughly like Figure~\ref{fig:below}.

\bigskip

\textbf{Second question}: what if one lets $K$ tend to infinity?

\bigskip

Since $\phi$ is conformal outside $\sq$ and $\phi(z)\sim z$ at infinity, the capacity of $\phi(\sq)$ is equal to that of $\sq$, and univalent function theory implies that the image $\phi(\sq)$ stays contained in some disk $B(0,R)$ independent of $K$. What else can be said?

There is a not-so-well-known similar case with a very simple solution: replace
$\sq$ with the unit disk $\D$, let $\mu$ be $0$ outside $\D$ and constant $\mu=a\in(0,1)$ within $\D$, corresponding to vertical ellipses with ratio $K=\frac{1+a}{1-a}$.
Then the straightening has the following form:
if $|z|\leq 1$ then $\phi(z)=z+a \bar{z}$, if $|z|\geq 1$ then $\phi(z)=z+a/z$.
The image $\phi(\D)$ of the unit disk is a horizontal ellipse, the lengthes of its semi major and semi minor axes being $1+a$ and $1-a$.
When $K\tend +\infty$, i.e.\ when $a\tend 1$, this ellipse tends to the segment $[-2,2]$ and $\phi$ converges uniformly. Its limit is $z+\bar{z}$ in $\overline{\D}$, which flattens $\D$ into $[-2,2]$, and $z+1/z$ on $\C\setminus \D$, which is a conformal mapping from the exterior of the unit circle to the complement of $[-2,2]$. 

\begin{figure}[htbp]
\begin{picture}(350,120)(0,0)
\put(0,10){\scalebox{0.5}{\includegraphics{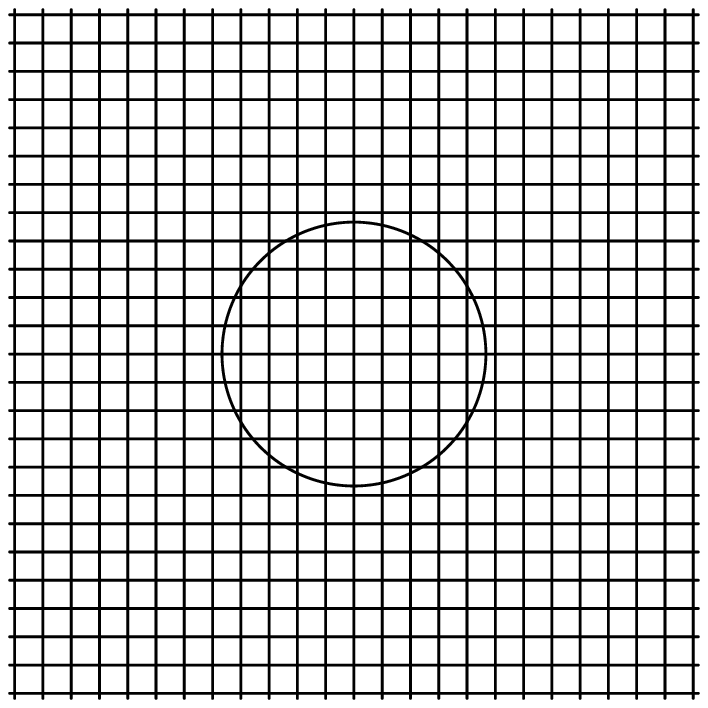}}}
\put(42,0){$K=1$}
\put(120,10){\scalebox{0.5}{\includegraphics{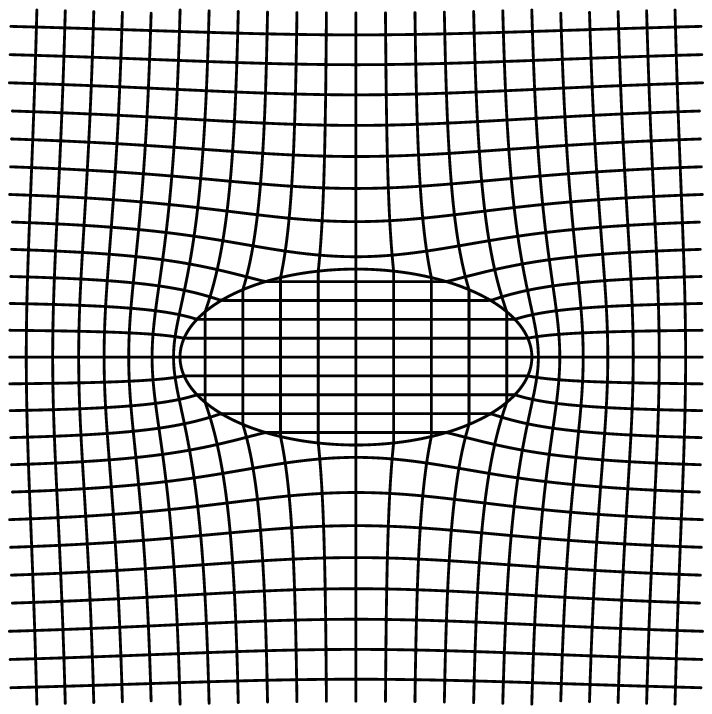}}}
\put(160,0){$K=2$}
\put(240,10){\scalebox{0.5}{\includegraphics{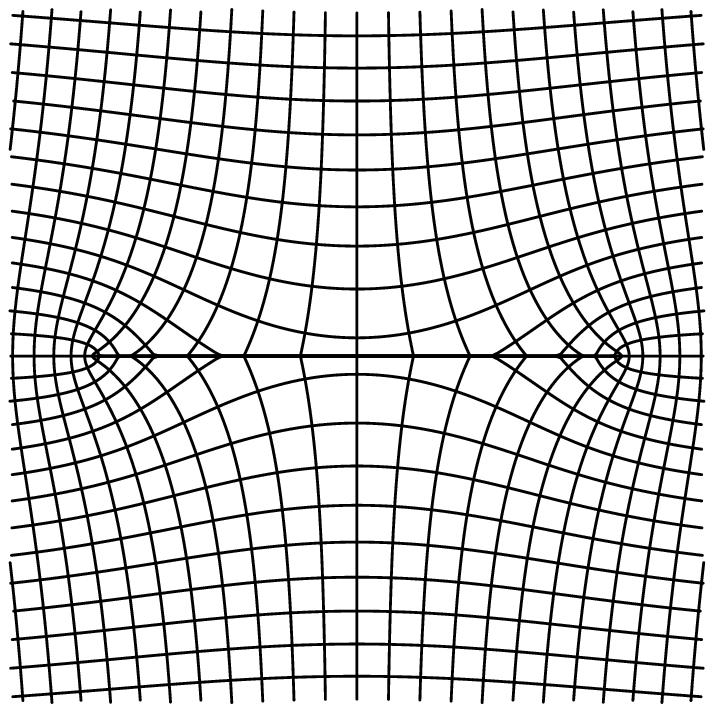}}}
\put(275,0){$K=+\infty$}
\end{picture}
\caption{Straightening a Beltrami form which is constant in the unit disk and zero outside.}
\end{figure}

But for the square, it cannot be so simple. For one thing, its right side segment accounts for $1/4$ of the harmonic measure\footnote{\emph{Harmonic Measure} has several meanings in mathematics, so let us precise the one we mean. Given a connected compact subset $K$ with more than one point and with connected complement, there is a unique $\rho>0$ and a unique conformal bijection $\phi$ from $\C\setminus K$ to $\C\setminus\rho\D$ such that $\phi(z)-z$ tends to $0$ at infinity. The potential associated to $K$ is the function $G(z)=\log |\phi(z)|$. It turns out that $G(z)=\int_{\partial K} \log|z-u| dm(u)$ for a unique measure $m$ with mass $1$ and support in $\partial K$, which is called the Harmonic measure with respect to $\infty$. It is also the unique non atomic measure with mass $1$ and support in $K$ (in fact $\partial K$) that minimizes the energy $E(m)=-\int_{K\times K} \log|z-w| dm(z) dm(w)$. The minimal energy is equal to $-\log \rho$, and is called the \emph{capacity}.} of $\partial \sq$ with respect to $\infty$. Therefore, its image by $\phi$ must account for $1/4$ of the harmonic measure of $\phi(\partial \sq)$ with respect to $\infty$. Thus it cannot tend to a point. Otherwise the capacity would have to tend to $+\infty$.

\subsection{First computer experiments}\label{subsec:exper}

To get a better idea of what is going on, the author, which will be referred to as `I' in the rest of this section, decided to make computer experiments, and took a Partial Differential Equation approach. Indeed, we can consider $\phi$ as a solution of $\Delta \phi = 0$, for a modified Laplacian $\Delta$.
The reader will see from the formulae below that the modified version still satisfies $\Re(\Delta\phi)=\Delta(\Re\phi)$ and $\Im(\Delta\phi)=\Delta(\Im\phi)$, so solving $\Delta \phi =0$ for the complex valued function $\phi$ amounts to independently solve the equation for its real and imaginary parts.

Far enough from $0$, $\phi(z)$ is close to $z$. So I took a big enough square, of side $R$, centered on $0$, and tried to numerically solve the following Dirichlet problem on a fine enough grid dividing this big square: $\Delta \phi=0$ with data $\phi(z)=z$ on the boundary of the big square for some \emph{modified discrete Laplacian} $\Delta$.

Usually, on a grid indexed by $(i,j)\in\Z^2$, given a function $(i,j)\mapsto u_{i,j}$, the standard discrete Laplacian has expression $\Delta u : (i,j) \mapsto u_{i-1,j}+u_{i+1,j}+u_{i,j-1}+u_{i,j+1}-4u_{i,j}$.
The modified discrete Laplacian is defined as follows. Let $z=z_{i,j}$ be the point on the grid corresponding to index $(i,j)$. This point $z_{i,j}$ can be outside $\sq$, inside $\sq$ or on one of its four sides $U$, $D$, $R$, $L$, or on a corner. Let us introduce the following local coordinates that locally solve the Beltrami equation near the point $(i,j)$: first we translate by $-z_{i,j}$ to put $z$ on the origin. Then we compose with a function that depends on where is $z$:
\begin{itemize}
  \item we compose with the identity if $z$ is outside $\sq$;
  \item with $(x+iy) \mapsto Kx+iy$ if $z$ is inside $\sq$;
  \item if $z$ is on the left side of the square $\sq$, we compose with $\on{id}$ on the left half-plane and $Kx+iy$ on the right half-plane;
  \item if $z$ is on the right side, we take them the other way round;
  \item if $z$ is on the upper side, we compose with $\on{id}$ in the upper half plane and $x+iy/K$ in the lower half plane;
  \item if $z$ is on the lower side, we take them the other way round;
  \item if $z$ is on a corner\ldots\ let us temporarily put corners aside.
\end{itemize}
In each case, let us call $l$, $r$, $u$, $d$ the length of the image, by this composition, of the segments from $z_{i,j}$ to the point $z_{i',j'}$ immediately on its left, on its right, above it, below it, respectively. The first idea was to define:
\[ (\Delta u)_{i,j}  =  \frac{ \ds\frac{u_{i+1,j}-u_{i,j}}{r} - \frac{u_{i,j}-u_{i-1,j}}{l}}{\ds\frac{r+l}{2}} + \frac{\ds\frac{u_{i,j+1}-u_{i,j}}{u} - \frac{u_{i,j}-u_{i,j-1}}{d}}{\ds\frac{u+d}{2}}
\]
Since the solutions of $\Delta \phi =0$ are the same if we multiply $\Delta$ by any function, I preferred to multiply this expression by $\frac{r+l}{2}\times\frac{u+d}{2}$ to get an expression that, you will notice, does not change when we compose the local coordinate with a $\C$-affine map. So let:
\begin{eqnarray*}
  (\Delta u)_{i,j} & = & \frac{u+d}{2r} (u_{i+1,j}-u_{i,j}) - \frac{u+d}{2l} (u_{i,j}-u_{i-1,j}) \\
   & & + \frac{l+r}{2u} (u_{i,j+1}-u_{i,j}) - \frac{l+r}{2d} (u_{i,j}-u_{i,j-1})
\end{eqnarray*}
Inside $\sq$ we get $(u_{i-1,j}+u_{i+1,j}-2u_{i,j})/K+(u_{i,j-1}+u_{i,j+1}-2u_{i,j})\times K$ and outside $\sq$ we just get the standard Laplacian. On its boundary it is a little bit more complicated. Now for corners: we take the standard Laplacian. It is not quite correct, but the heuristics is that getting the formula wrong on a bounded number of points will not influence too much the result, as long as we take something reasonable: here $(\Delta u)_{i,j}$ vanishes if and only if the central value is a barycenter of the neighborhing values, (the weights of the barycenter are allowed to depend on $i,j$).

\begin{figure}[htbp]
\begin{picture}(350,175)
\put(0,0){\rotatebox{0}{\scalebox{0.4}{\includegraphics{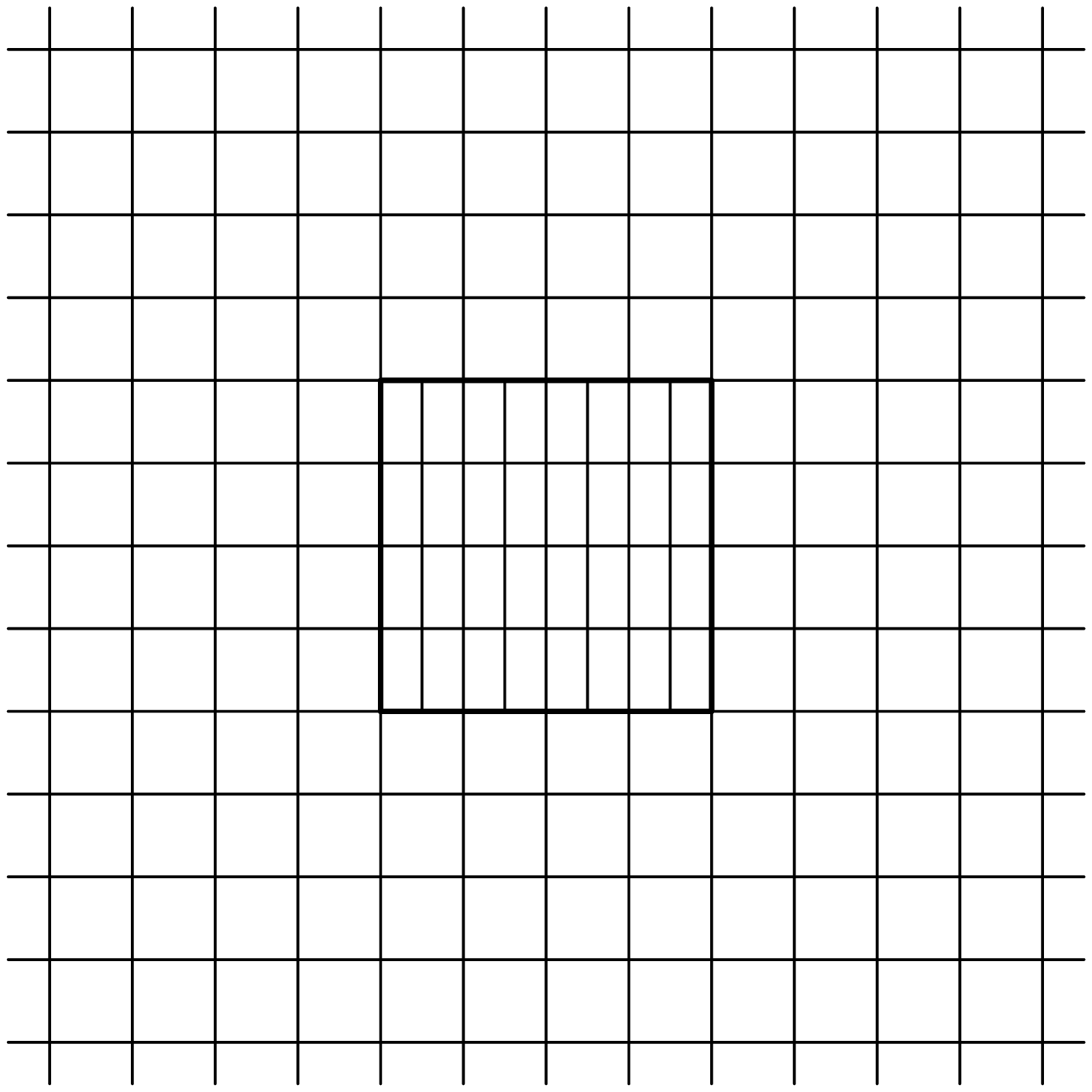}}}}
\put(175,5){\rotatebox{90}{\scalebox{0.31}{\includegraphics{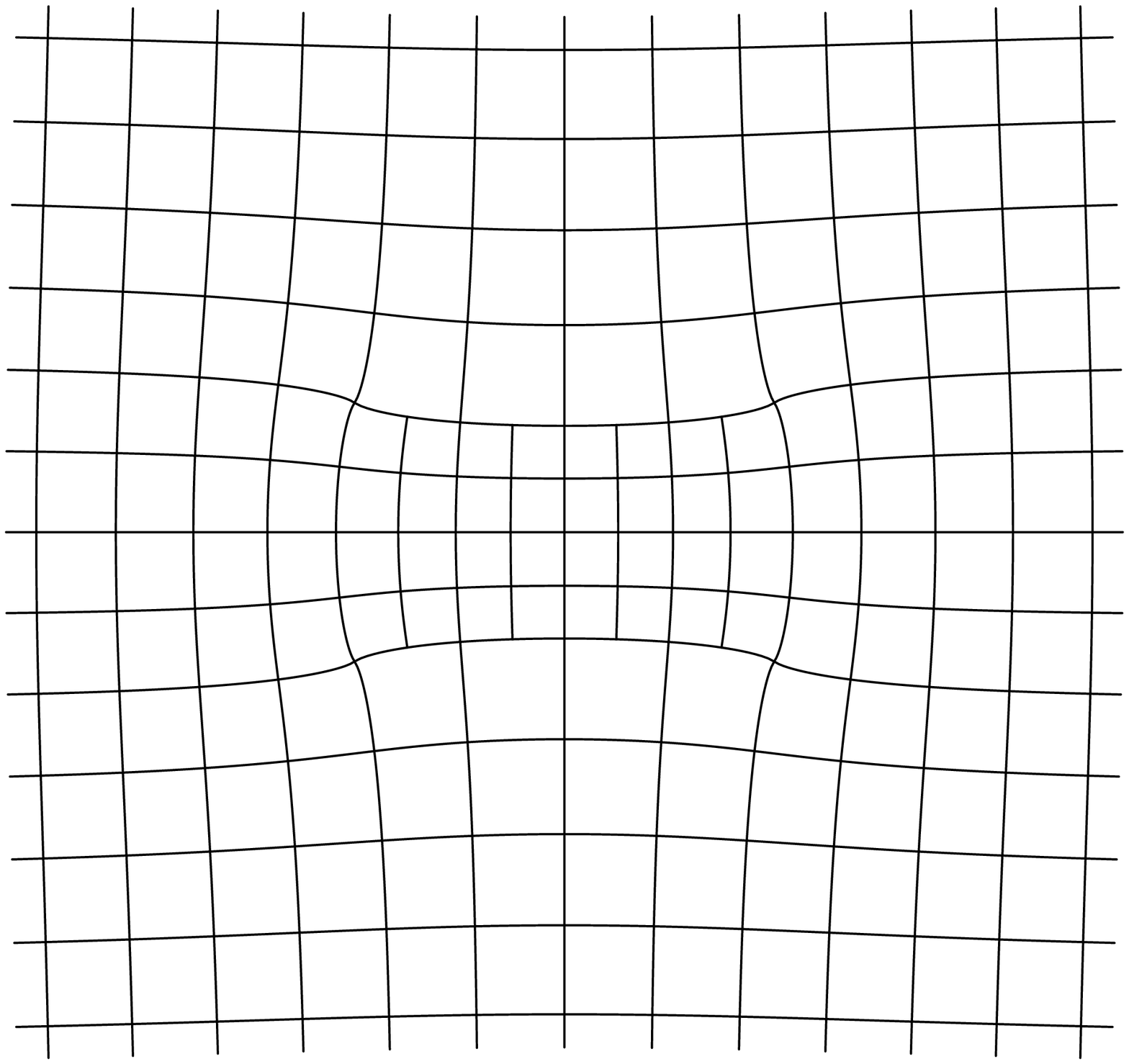}}}}
\end{picture}
\caption{The result of the computer experiment. Left: a grid, which within the square $\sq$ has been divided into rectangles of ratio $2$. Right: the image of the left grid by $\phi$. The rectangles have become conformal squares. \footnotesize{Note: The grid $G$ on which the PDE was solved to produce this picture is finer than the one which is drawn. The spacing of $G$ is $1/64$ times the length of the side of the square supporting the Beltrami form. And $G$ extends beyond the image: it extends in the 4 directions to a distance of 10 times the square's side.}}
\end{figure}

Let me now describe the numerical method I used to solve the discrete equation $\Delta \phi=0$ with the computer.\footnote{Acknowledging I am not an expert in numerical solutions of PDEs, there is probably room for improvement there.} I used the \emph{Jacobi Relaxation Method}: recall that $\Delta u$ vanishes at the vertex $(i,j)$ if and only if $u_{i,j}$ is a barycenter of the value of the neighboring vertices (with weights that depend only on $(i,j)$). The idea is then to replace $u_{i,j}$ by this barycenter. 
More precisely, define inductively a sequence $u^{(n)}$ of functions $(i,j)\mapsto u^{(n)}_{i,j}$ as follows: begin with $u^{(0)}_{i,j}=z_{i,j}$, which agrees with imposed boundary values. Now given $u^{(n)}$, if $(i,j)$ is not a boundary vertex, let $u^{(n+1)}_{i,j}$ be the value of the aforementioned barycenter computed from the values of $u^{(n)}$ at the vertices neighboring $(i,j)$. If $(i,j)$ is a boundary vertex, define $u^{(n+1)}_{i,j}=u^{(n)}_{i,j}$. The reader will notice that since we change the value everywhere (but on the boundary), the Laplacian of $u^{(n+1)}$ is still not $0$. However, the values of $u^{(n)}$ are supposed to converge to the solution of $\Delta u = 0$.

\begin{remark} I made no effort at trying to prove the the convergence of this scheme and the correctness of its limit, nor at estimating the speed of convergence.
\end{remark}

As this was a little bit slow, I first solved the equation on a coarse grid, and then used the (near) solution $u^{(N)}$, for some big $N$, to define a starting point $v^{(0)}$ for solving the equation on a finer grid (with the grid unit divided by $2$, the value of $v^{(0)}$ on the newly introduced points being defined by averaging on its neighbors), and went on refining the grid further several times. This resulted in a great acceleration of the convergence.

\begin{figure}[htbp]
\begin{picture}(348,306)
\put(0,153){\rotatebox{90}{\scalebox{0.3}{\includegraphics{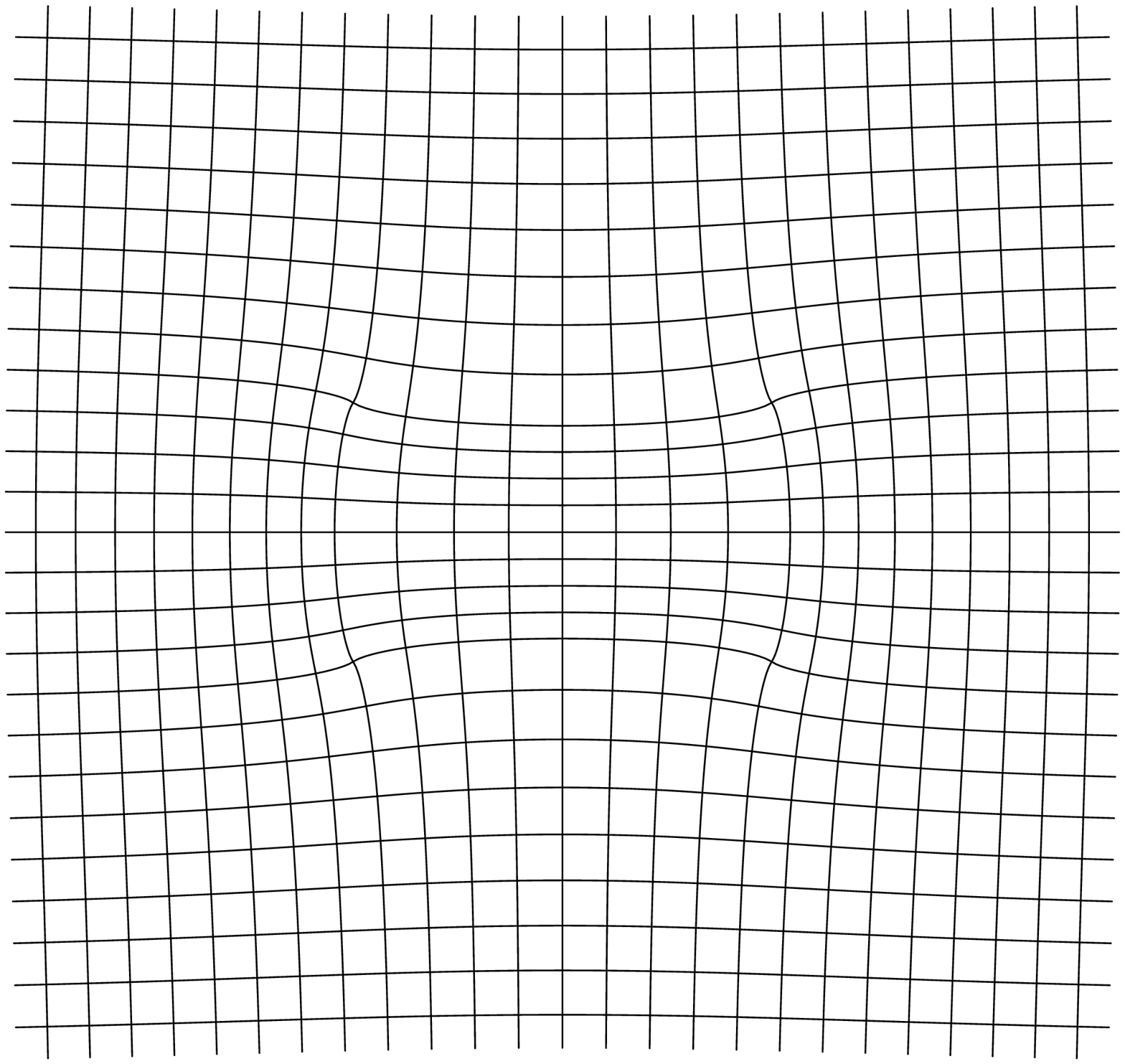}}}}
\put(174,153){\rotatebox{90}{\scalebox{0.3}{\includegraphics{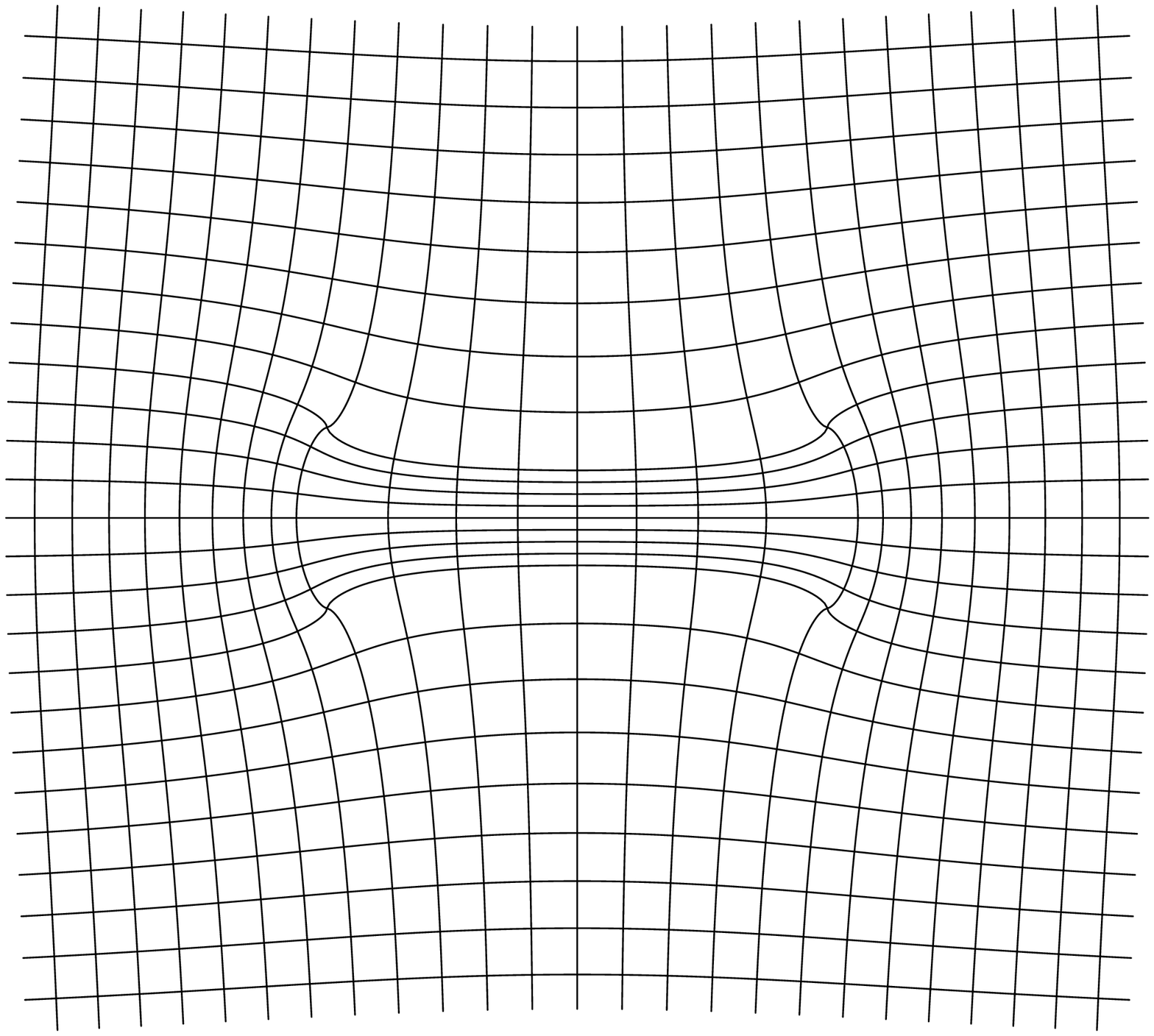}}}}
\put(0,0){\rotatebox{90}{\scalebox{0.3}{\includegraphics{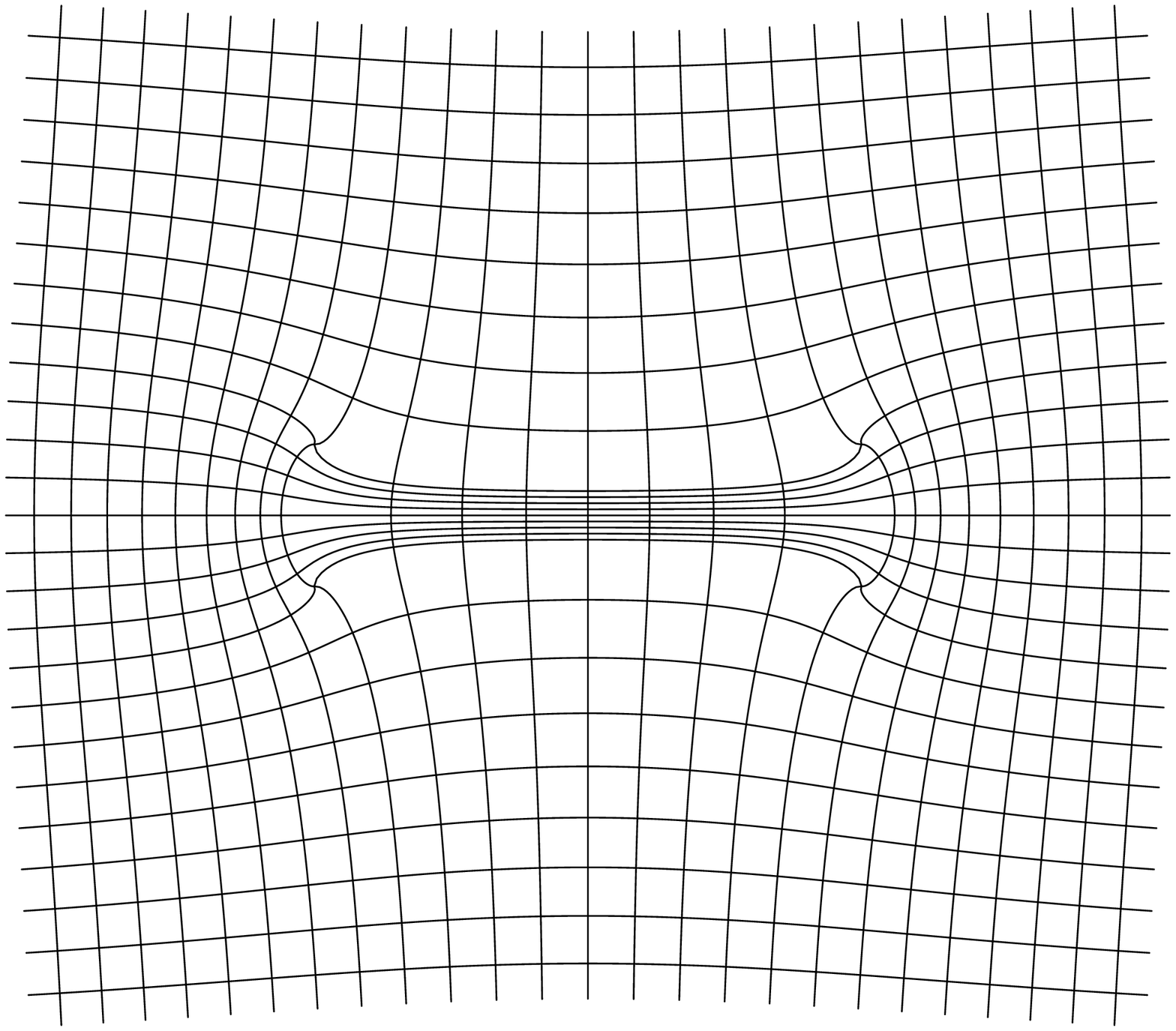}}}}
\put(174,0){\rotatebox{90}{\scalebox{0.3}{\includegraphics{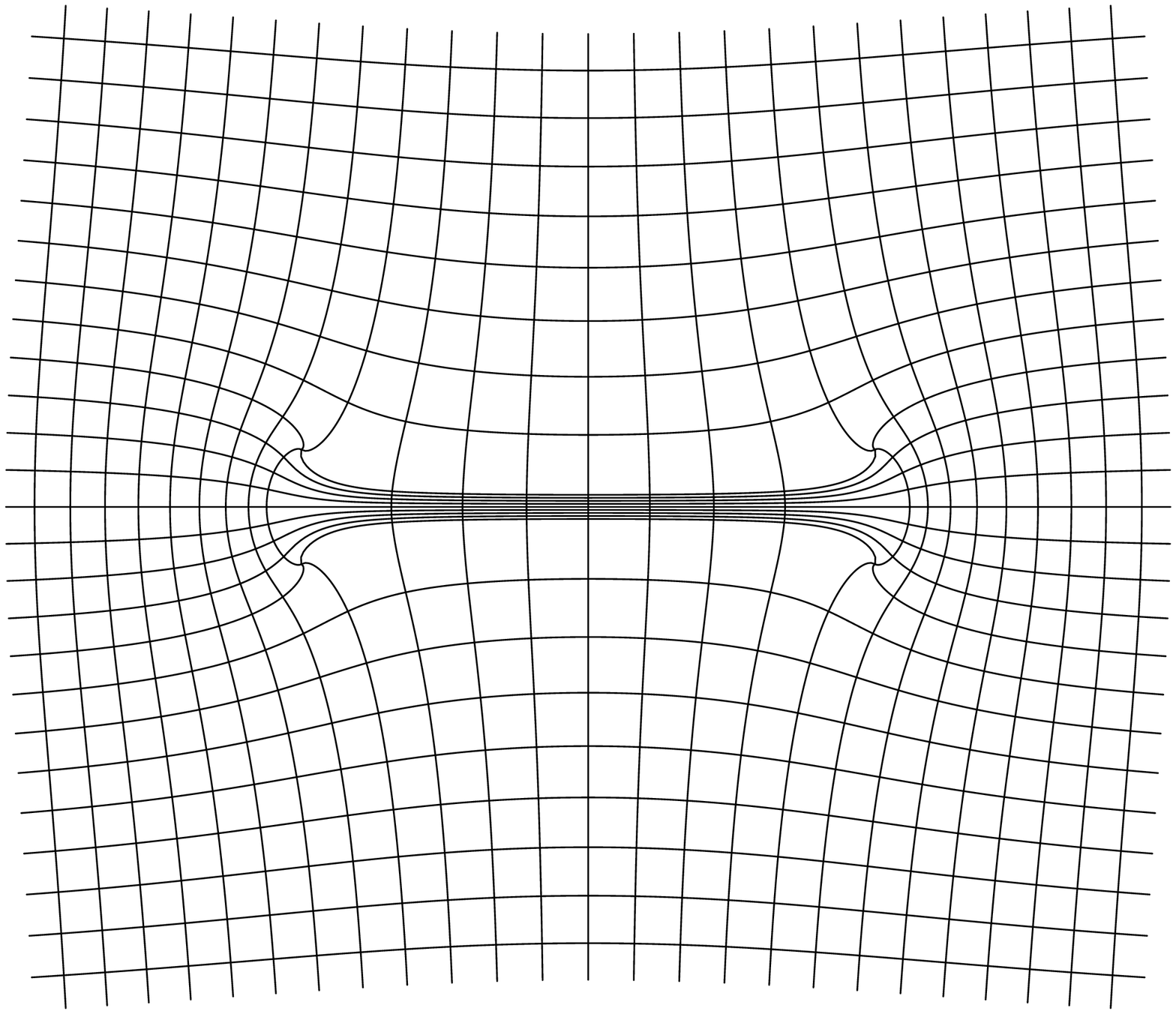}}}}
\end{picture}
\caption{The image of the square, as approximated by the program, for $K=2$, $5$, $10$, $20$.}
\label{fig:a}
\end{figure}

Let us list a few defects of the method:
\begin{itemize}
\item The definition of the modified Laplacian makes use of the fact that the side of the square and the major axis of the ellipses are parallel to the coordinate axes. It is not clear how to adapt it to more general ellipse fields.
\item It is computer intensive.
\item It may have a too high algorithmic complexity (the time it takes to reach a given accuracy $\epsilon$).
\item The grid is uniform, whereas a non uniform mesh would be more efficient.
\item When $K$ tends to infinity one needs at the same time a finer grid and a bigger $N$. The Laplacian approach does not seem well suited to big values of $K$.
\end{itemize}
The reader will find on figure~\ref{fig:a} what the same program yields when we increase $K$. It was hard to get something convincing in a reasonable amount of time for $K=50$. 
A few hints for improvement could be:
\begin{itemize}
\item One could have used the Conjugate Gradient Method (see \cite{Jo}) instead of Jacobi Relaxation.
\item One may use a better approximation than $\phi(z)=z$ on the boundary of the big square, allowing to take a smaller one.
\item Replacing, within the square, the $n\times n$ grid by an $n\times m$ grid, with $n/m\approx K$.
\item Using a Finite Elements Method.
\end{itemize}
All this discussion was about solving $\Delta \phi=0$. Amongst other approaches, there are:
\begin{itemize}
\item Circle Packings (see \cite{He} and \cite{Brw} for instance). However, it probably needs a very dense mesh too when $K$ tends to infinity.
\item Fourier series and the Hilbert-Beurling singular integral operator. See for instance \cite{Dar}, \cite{GK}. One could also periodize the problem (using the big square as a fundamental domain) and use standard 2D Fourier series on the torus.
\item Discrete Riemann surfaces. See \cite{Mer}.
\end{itemize}
\par\medskip
But I did not push my experiments further. In any case, big values of $K$ seemed to require a lot of work and computer time. The behavior of $\phi(\sq)$ when $K\tend +\infty$ remained a mystery for me. I was short of theoretical arguments, and there was no hope of an explicit formula there, unlike the case of the disk\ldots\ Or so I thought.

\begin{figure}[htbp]
\scalebox{0.65}{\includegraphics{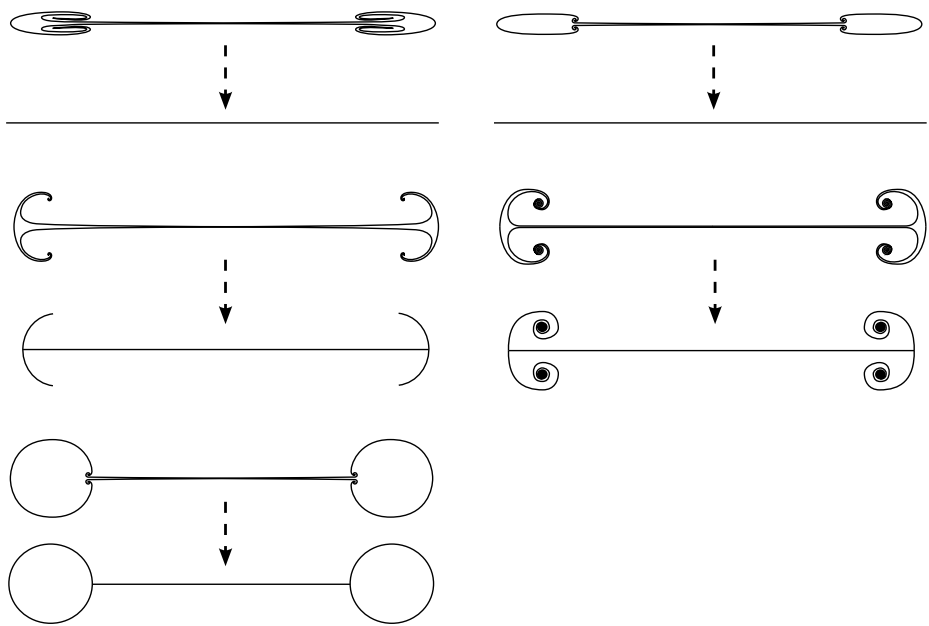}}
\caption{Five guesses of the possible limit shape, if there is one, together with a possible way to tend to it.}
\label{fig:c}
\end{figure}

\subsection{A few guesses}\label{subsec:guess}

Imagining that $\phi(\sq)$ has a limit when $K\tend\infty$, what could it be? It should be noted that the spiral, which was quite slowly turning for $K=2$, does it faster and faster. But does its diameter tend to $0$? 
Figure~\ref{fig:c} gives a few possibilities. We will see later in the article that the reality is even more\ldots\ interesting.

\section{An explicit formula}\label{sec:sol}

It turns out that there is an (almost completely) explicit formula for $\phi^{-1}$. 
It allowed the author to draw very precisely on the computer the image of the square, and to push the value of $K$ up to $10^{50}$ and beyond\ldots\ Thanks to these pictures, it was possible to guess the limit of $\phi(\sq)$, at least qualitatively. The formula can be proved to have a limit that allows to give a guess of the actual limit of $\phi(\sq)$. We will see that the theoretical framework used to discover the formula gives an interesting interpretation of these limits.

\begin{remark} What we will explain here generalizes to the straightening of Beltrami forms defined on $\C$, piecewise constant, with polygonal pieces, and finitely many of them. This will be covered in a forthcoming article.
\end{remark}

\begin{figure}[htbp]
\begin{picture}(315,170)
\put(0,0){\includegraphics{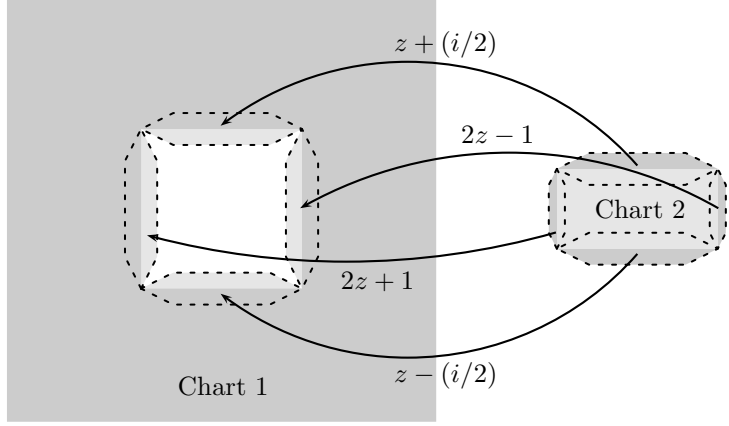}}%
\put(69,15){Chart 1}
\put(225,82){Chart 2}
\put(155,27){\psarc{->}(0,0){4}{41}{128}}
\put(150,145){$z+(i/2)$}
\put(155,143){\psarc{<-}(0,0){4}{232}{319}}
\put(150,20){$z-(i/2)$}
\put(193,-50){\psarc{->}(0,0){5.5}{60}{120}}
\put(175,110){$2z-1$}
\put(132,350){\psarc{<-}(0,0){10}{254.8}{286}}
\put(130,55){$2z+1$}
\end{picture}
\caption{The affine surface behind the formula. Case $K=2$. The arrows are $\C$-affine maps (similitudes).}
\label{fig:d}
\end{figure}

\titre{Defining an affine surface} Consider the following Riemann surface: one chart is given by the exterior of $\sq$, union $\partial \sq$, minus the four corners. The second chart is a rectangle given by the image of $\sq$ minus its corners by $x+iy \mapsto x+i(y/K)$. Now glue the four sides of the first chart to the four sides of the second chart by four affine maps, as illustrated on figure~\ref{fig:d}.\footnote{Normally, charts should be open. So what is meant is that for all sufficiently small open sets containing our two charts, gluing them with the previous affine maps (more precisely by the unique $\C$-affine extensions), we get Riemann surfaces that are all ``the same'', i.e.\ isomorphic by the obvious map.} Not only we get a Riemann surface, but the transition maps are $\C$-affine. So we get an \emph{affine surface}\footnote{An affine surface is a topological surface with an atlas whose transition maps are all locally $\C$-affine maps, i.e.\ of the form $z\mapsto az+b$ (being locally $\C$-affine on an open set, the transition map is affine on each connected component of its domain of definition). An affine surface is in particular a Riemann surface by the same atlas.} $\cal S$. 

\titre{Completing the Riemann surface} As a topological surface $\cal S$ is homeomorphic to the sphere minus five points. They correspond to $\infty$ and the four corners of the square. Let us prove that we can extend the Riemann surface $\cal S$ to these five points,\footnote{Consider a topological surface $S$ and a Riemann surface structure on an open subset $U$. Assume there is a point $x\in S$ that does not belong to $U$ but has a neighborhood included in $U\cup\{x\}$. It is \emph{not always possible} to extend the Riemann surface structure to $U\cup\{x\}$.} by providing adequate charts near them. The easiest one to add is $\infty$: we introduce the chart $B(0,\epsilon)$ and glue it to Chart~$1$ with the map $z\mapsto 1/z$. 
To add a corner, remark that a neighborhood of it is isomorphic, as an affine surface, to a neighborhood of $0$ of the following affine surface: $\C$ with a straight slit radiating away from $0$, the two sides of the slits being glued together by a map of the form $z\mapsto Kz$. We have already seen in section~\ref{subsec:pb} how to conformally map this to $\C^*$: a branch of the map $z\mapsto z^\alpha$ will work, with 
\begin{eqnarray*}
\alpha &=& 1/(1-\log(K)/2i\pi)\text{ for the upper right and lower left corners,}\\
 \alpha &=& 1/(1+\log(K)/2i\pi)\text{ for the other two.}
\end{eqnarray*}
This gives a conformal chart in which we can add the origin, i.e.\ the corner. Therefore we have proved:
\begin{lemma}\label{lem:loccorner}
  Each corner in $\cal S$ has a neighborhood that can be conformally mapped to a neighborhood of the origin in $\C$, on which any branch of $z\mapsto z^{1/\alpha}$ is an affine chart (for the value of $\alpha$ defined above).
\end{lemma}

\titre{Uniformizing to the Riemann sphere} We have thus extended the Riemann surface $\cal S$ into a compact Riemann surface $\wt{\cal S}$, with five more points. Since $\wt{\cal S}$ is homeomorphic to the sphere, it is conformally equivalent to the Riemann sphere $\S$.
Let $\Phi:\wt{\cal S}\to \S$ be the unique isomorphism normalized by the following condition: 
\begin{equation}\label{eq:normali}
  \phi_1(z)=z+0+o(1)\text{ as }z\tend \infty
\end{equation}
where $\phi_1$ is the expression\footnote{i.e.\ $\phi_1=\Phi\circ\on{chart}_1^{-1}$ where $\on{chart}_1$ goes from the corresponding open subset of $\cal S$ to $\C\setminus\sq$} of $\Phi$ in Chart~$1$.
The subset $\cal S$ of $\wt{\cal S}$ is mapped by $\Phi$ to the complex plane minus four points. By the symmetries of the surface, the chosen normalization, and uniqueness of $\Phi$, these points must be of the form $z_1=x+iy$, $z_2=-x+iy$, $z_3=-x-iy$, $z_4=x-iy$, for some $x>0$ and $y>0$.

\titre{Connection with the straightening} We can use the affine surface $\cal S$, and the uniformization $\Phi$ of its completion, to recover the straightening $\phi$ of the Beltrami form that is our object of interest. Let $\Psi:\C\setminus\{1+i,-1+i,-1-i,1-i\}\to\cal S$ be defined as follows: $\Psi$ maps a point $z$ outside the square $\sq$ to the point of $\cal S$ which has coordinate $z$ in Chart~1. It maps a point $z=x+iy$ within the square to the point of coordinates $x+i(y/K)$ in Chart~2. Then $\Phi\circ\Psi$ is a homeomorphism from $\C\to\C$, that satisfies $\Phi\circ\Psi(z)=z+0+o(1)$ as $z\tend\infty$. It is a $C^1$ real diffeomorphism on the complement of the boundary of the square, that straightens the Beltrami form $\mu$. By quasiconformal erasability of lines, it is the straightening of $\mu$:
\begin{equation}\phi=\Phi\circ\Psi.
\end{equation}

\titre{Equivalence of an affine surface structure compatible with a Riemann surface structure, and the distortion derivative}
The affine surface structure on $\cal S$, i.e.\ the notion of affine chart, can be transported by $\Phi$ to an affine surface structure on $\C\setminus\{z_1,z_2,z_3,z_4\}$. 
Indeed, it is given by some atlas on $\cal S$, and by composing charts with $\Phi^{-1}$, we get an atlas on $\C\setminus\{z_1,z_2,z_3,z_4\}$.
Now consider any affine chart $(U,f)$ on $\C\setminus\{z_1,z_2,z_3,z_4\}$, where $U$ is an open subset of $\C\setminus\{z_1,z_2,z_3,z_4\}$ and $f : U\to \C$. The map $f$ is analytic, because $\Phi$ is. Since it is also injective, $f'$ does not vanish. Since the transition maps between affine charts are affine, the complex number $\eta(z)=f''(z)/f'(z)$, called the \emph{non-linearity} or the \emph{distortion derivative} of $f$, does not depend on the choice of the chart $(U,f)$ near $z$. Indeed for any two affine charts $f_1 : U_1 \to \C$ and $f_2 : U_2 \to \C$, and any $z\in U_1\cap U_2$, there exists $a,b\in\C$ with $a\neq 0$ such that $f_2=a f_1+b$ holds near $z$. Reciprocally, from $\eta(z)$ it is possible, at least locally, to recover the affine charts: $f=\int \exp \left(\int \eta\right)$, the two integration constants accounting for $a$ and $b$. In other words the function $\eta$ characterizes the affine structure.

We therefore have a well-defined holomorphic function $\eta: \C\setminus\{z_1,z_2,z_3,z_4\} \to\C$, satisfying $\eta(z)=f''(z)/f'(z)$ for all affine charts $f$ near $z$, and completely characterizing the affine structure. By equation~\eqref{eq:normali}, since $f=\phi_1^{-1}$ is an affine chart, we get that $\eta(z) \tend 0$ as $z\tend \infty$.

\titre{Behavior of the distortion derivative under a change of variable}
To study the singularity at $z_i$ of $\eta$, we will make a convenient change of variable: by Lemma~\ref{lem:loccorner} there exists a conformal map $\psi$ from a neighborhood of $z_i$ to a neighborhood of $0$ in $\C$ such that for any branch $g$ of $z\mapsto z^{1/\alpha}$, the map $f=g\circ\psi$ is an affine chart of $\C\setminus\{z_1,z_2,z_3,z_4\}$ near $z_i$; the number $\alpha$ equals $1/(1-\log(K)/2i\pi)$ for $z_1$ and $z_3$, and $1/(1+\log(K)/2i\pi)$ for $z_2$ and $z_4$.
The quantity $g''(z)/g'(z)$ that characterizes the affine structure in the new conformal chart is equal to $\beta/z$ where $\beta=\alpha^{-1}-1$. At $z_1$ and $z_3$, $\beta = -\log(K)/2i\pi$. At $z_2$ and $z_4$, $\beta = \log(K)/2i\pi$.
Since $f=g\circ\psi$, we have the composition formula:
\begin{equation}\label{eq:cocy}
 \frac{f''}{f'} = \psi' \times \frac{g''}{g'}\circ \psi + \frac{\psi''}{\psi'}.
\end{equation}
Now $\psi'$ does not vanish at the origin, and it follows from this and from $g''/g'=\beta/z$ that
the map $\eta=f''/f'$ has a simple pole at $z_i$ with residue\footnote{In fact the equality of the residues of $f''/f'$ and $g''/g'$ still holds if the pole has higher order or is an essential singularity, as can be seen by multiplying equation~\eqref{eq:cocy} by $dz$ and integrating both sides on a small loop around $z_i$.} $\beta$.

\titre{The distortion derivative is completely determined} Indeed $\eta$ is a holomorphic function over $\C\setminus\{z_1,z_2,z_3,z_4\}$ that tends to $0$ at $\infty$, and has simple poles poles at $z_1$, $z_2$, $z_3$, $z_4$ with residue $-\log(K)/2i\pi$ at $z_1$ and $z_3$ and  $\log(K)/2i\pi$ at $z_2$ and $z_4$. Such a function is necessarily rational and:
\begin{equation} \eta(z) = -\frac{\log K}{2i\pi} \left( \frac{1}{z-z_1} - \frac{1}{z-z_2} + \frac{1}{z-z_3} - \frac{1}{z-z_4} \right).
\end{equation}
So we know $\eta$ explicitly (and thus the affine charts $f$, and from this the inverse $\phi^{-1}$ of the straightening) \emph{provided} we can determine the value of $z_i$. This can be done numerically, far much faster than for a Partial Differential Equation, as we will see in section~\ref{sec:result}.

\titre{Local expression of the affine coordinates} Explicitly solving the equation $f''/f'=\eta$ yields: 
\begin{eqnarray*}
\log f'  & = &  \frac{\log K}{2i\pi} \Big(\log(z-z_4)-\log(z-z_3)+\log(z-z_2)-\log(z-z_1) \Big)+\text{const}
\\ & = & \frac{\log K}{2i\pi} \log\left(\frac{(z-z_2)(z-z_4)}{(z-z_1)(z-z_3)}\right) + \text{const}
\end{eqnarray*}
where we can guess there will be branch problems.
At least locally the charts $(U,f)$ on $\C\setminus\{z_1,z_2,z_3,z_4\}$ have the expression
\begin{equation}
  f(z) = a \int \left( \frac{(z-z_2)(z-z_4)}{(z-z_1)(z-z_3)} \right)^{\frac{\log K}{2i\pi}}dz +b
\end{equation}
for some $a\in\C^*$ and $b\in\C$. Determining the two affine charts corresponding to Chart~1 and Chart~2 thus requires to determine the value of $a$ and $b$ for each.

\titre{Going global} The statement above is local. To extend it a little bit, we must be careful about monodromy questions. For one thing, the integrand 
\[ m(z)=\left( \frac{(z-z_2)(z-z_4)}{(z-z_1)(z-z_3)} \right)^{\frac{\log K}{2i\pi}} \]
is a multivaluated expression: the set of possible values at a given $z$ is of the form $\setof{K^n v}{n\in\Z}$ for some $v\in\C^*$. Also, any affine chart $f$ with a connected domain of definition extends to any bigger and simply connected open subset of $\C\setminus\{z_1,z_2,z_3,z_4\}$ into a unique solution of $f''/f'=\eta$, but the extension is not necessarily injective. 

\titre{Extending an affine chart along a path}
Consider a path that is contained in $\C\setminus\{z_1,z_2,z_3,z_4\}$, ends included. Let $z_0$ be its starting point. Consider the initial values $f'(z_0)$ and $f(z_0)$. We can choose any determination $m_0$ of $m(z_0)$ and let $b=f(z_0)$ and $a=f'(z_0)/m_0$. Then the integrand $m(z)$ can be followed continuously along the path and its integral along the path is well defined, which gives a certain value for the expression $a\int m(z)dz +b$. Any homotopic path will yield the same value. Let us call this the continuation of the affine chart $f$ along the path. It depends only on $f(z_0)$, $f'(z_0)$ and the homotopy class of the path.\footnote{Thus we have a well defined analytic map on the universal cover of $\C\setminus\{z_1,z_2,z_3,z_4\}$. It is called the \emph{developing map} of $f$.}
If the path is completely contained in the domain of definition $\on{Def}(f)$ of an affine chart $f$ (or the domain of definition of an analytic continuation $f$ of an affine chart), then we recover the value of $f$ at the end of the path. This holds even if this domain is not simply connected: the mere existence of $f$ implies that there is no monodromy within $\on{Def}(f)$. In particular, there exists a branch of $m$ defined on $\on{Def}(f)$.
If the path is not completely contained by just homotopic to a path completely contained in $\on{Def}(f)$, we recover the same value.

\titre{Expression of Chart~1} Identifying $\cal S$ with $\C\setminus\{z_1,\cdots,z_4\}$ via $\Phi$, Chart~1 can be considered as defined on $U_1=\phi(\C\setminus\sq)$ and equal there to $\phi^{-1}$, since $\Psi=\on{id}$ on the complement of $\sq$. Let us choose the branch of $z\mapsto m(z)=\left( \frac{(z-z_2)(z-z_4)}{(z-z_1)(z-z_3)} \right)^{\frac{\log K}{2i\pi}}$ on $U_1$ that tends to $1$ at $\infty$. According to the above analysis, $\phi^{-1} = b + a\times($ a primitive of $m$ on $U_1)$. Since $\phi^{-1}{}'(z)\tend 1$ at $\infty$, we must have $a=1$. Note that $m(z)-1$ is necessarily integrable on pathes within $U$ that tend to $\infty$, because $\phi^{-1}(z)-z$ has a limit at $\infty$ and $m(z)=(\phi^{-1})'$. (In fact, $m(z)=1+\cal O(z^{-2})$ at $\infty$, as can be shown by a direct Taylor expansion, using the fact that $z_2+z_4=0$ and $z_1+z_3=0$.)
Since $\phi^{-1}(z)-z$ tends to $0$, we get:

\begin{proposition}\label{prop:expr1} The straightening $\phi$ of the Beltrami form defined in section~\ref{subsec:pb} has the following property outside of the square $\sq$: for all $z\in U_1=\phi(\C\setminus\sq)$, we have
  \[\phi^{-1}(z) = z+\int_\infty^z \left( \left(\frac{(u-z_2)(u-z_4)}{(u-z_1)(u-z_3)} \right)^{\frac{\log K}{2i\pi}} -1\right) du,\]
  on any path from $\infty$ to $z$ that is contained in $U_1$ (or just homotopic within $\S\setminus\{z_1,z_2,z_3,z_4\}$ to a path contained in $U_1$ with the same endpoints), for the branch of the expression $\left( \frac{(z-z_2)(z-z_4)}{(z-z_1)(z-z_3)} \right)^{\frac{\log K}{2i\pi}}$ defined on the path and equal to $1$ at $\infty$.
\end{proposition}

\begin{remark} There is a difficulty in using this statement: it does not tell us how to determine whether or not a given $z$ belongs to $U_1$. Moreover, remember that $U_1$ spirals around the points $z_i$. Therefore, given some $z\in U_1$, it is not obvious how to determine if a path going to $\infty$ is homotopic in the complement of the $z_i$'s to one staying in $U_1$.
\end{remark}

\smallskip

\titre{Expression of Chart~2} Using Chart~2 translated by $i(1-K^{-1})$, we can extend Chart~1 to $U=\phi((\C\setminus\partial \sq)\cup(-1+i,1+i))$. It corresponds to gluing Chart~2 to Chart~1 along the upper side of the square.
According to the analysis above, the branch of $m$ we have chosen on $U_1$ extends to $U$, since the affine chart does.
\begin{figure}[htbp]
\includegraphics[height=5cm]{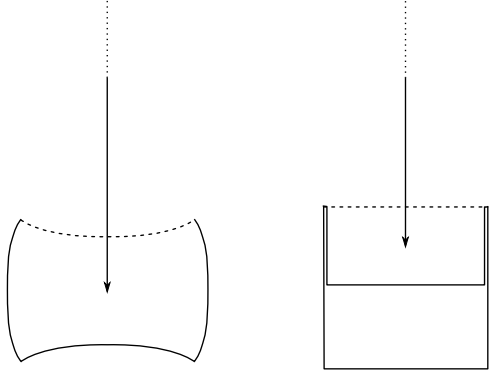}
\caption{Extending Chart~1 to $U=\phi((\C\setminus\partial \sq)\cup(-1+i,1+i))$.}\label{fig:guess3}
\end{figure}
Let $\theta_1$ be the unique argument of $z_1$ that belongs to $(0,\pi/2)$.
Notice that at $z=0$, the quantity $\frac{(z-z_2)(z-z_4)}{(z-z_1)(z-z_3)}$ takes the value $e^{-i4\theta_1}$.
The path going from $\infty$ down to $0$ within the imaginary axis is contained in $U$. Following the logarithm of $\frac{(z-z_2)(z-z_4)}{(z-z_1)(z-z_3)}$ along the path, starting with the determination $0$ of $\log 1$, leads to the value $-i 4\theta_1$.
From this, it follows that $m(0)=\exp(-i4\theta_1 \log(K)/2i\pi) = \exp(-2\theta_1\log(K)/\pi)$.

\begin{proposition}\label{prop:expr2}
 For all $z\in U_2=\phi(\sq)$, we have
  \[\phi^{-1}(z) = A\Bigg(\int_0^z \bigg(
    \frac{(u-z_2)(u-z_4)}{(u-z_1)(u-z_3)} 
    \bigg)^{\frac{\log K}{2i\pi}} du\Bigg),
  \]
  where $A(x+iy)=x+iKy$, and the path of integration between $0$ and $z$ is contained in $U_2$ (or just homotopic within $\S\setminus\{z_1,z_2,z_3,z_4\}$ to a path contained in $U_2$ with the same endpoints), for the branch of the expression $\left( \frac{(u-z_2)(u-z_4)}{(u-z_1)(u-z_3)} \right)^{\frac{\log K}{2i\pi}}$ defined on the path and taking at $z=0$ the value $\exp(-2\theta_1\log(K)/\pi)$, where $\theta_1$ is the unique argument of $z_1$ that belongs to $(0,\pi/2)$.
\end{proposition}

\begin{remark} We have the same difficulties as in the previous proposition.
\end{remark}

In the next section, we will explain how to practically use these formulae and how to determine the value of $z_1$ and the set $\phi(\sq)$.

\section{Practical aspects}\label{sec:prac}

\subsection{Integrating towards a singularity}

Consider the continuation of an affine chart along a path, as defined in the previous section.
The following lemma tells that under certain conditions, the path can tend to one of the corners $z_i$.
\begin{figure}[htbp]
\includegraphics[height=3cm]{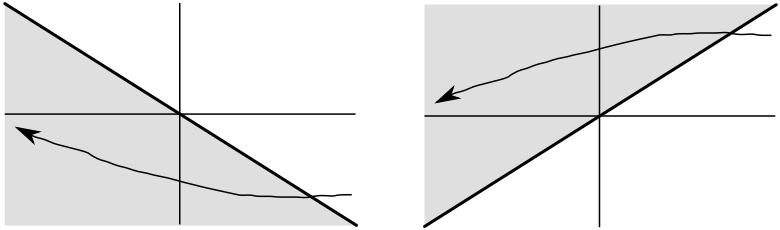}
\qquad
\rotatebox[origin=c]{180}{\reflectbox{\includegraphics[height=3cm]{img/hp.eps}}}
\caption{The half plane $\Re( z) -\tau \Im(z) < 0$. Left: when $\tau<0$. Right: when $\tau>0$. A lift $L(t)$ of the path $\gamma(t)$, satisfying $\exp(L(t))=\gamma(t)-z_i$ is drawn.
The real part of the lift of the path must tend to $-\infty$ and its distance to the line must tend to $+\infty$.}
\label{fig:hp}
\end{figure}
\begin{lemma}\label{lem:sing}
 Consider an initial point point $z_0 \in \C\setminus\{z_1,z_2,z_3,z_4\}$, and a path $\gamma:[0,1]\to\C$ from $z_0$ to one of the corners $z_i$, that is contained in $\C\setminus\{z_1,z_2,z_3,z_4\}$ except at its arrival.
Let $L:[0,1)\to\C$ be a continuous lift such that $\exp(L(t))=\gamma(t)-z_i$.
Let $\tau = \log(K)/2\pi$ if $i=1$ or $3$, and $\tau=-\log (K)/2\pi$ if $i=2$ or $4$.
Then a continuation of an affine chart along $\gamma$ will converge at $t=1$ if and only if
\[\Re(L(t)) - \tau\Im(L(t)) \underset{t\to 1^-}{\tend} -\infty\]
(see Figure~\ref{fig:hp}).
\par
Consider two paths starting from $z_0$, ending at $z_i$, satisfying the above condition on their lifts, and homotopic within $\C\setminus\{z_1,z_2,z_3,z_4\}$. The continuation of an affine chart along these paths will converge and tend to the same value at $t=1$.
\end{lemma}
\begin{proof} We may post-compose the affine chart with an affine map.
\par
First claim: Let $\beta=i\tau$.
We saw that there is a change of variable $z=\psi(w)$ defined near $z_i$, sending $z_i$ to $0$ and such that, locally, the affine charts $g=f\circ\psi$ satisfy $g''(w)/g'(w) = \beta/w$ (and thus are are branches of $w\mapsto a w^{\beta+1} +b$). This affine structure can then be lifted by the exponential: let $w=\exp(u)$ and $h = g \circ \exp$. Then $h''/h'=\beta+1$. This equation has global solutions: let us choose $h(z)=\exp\big((\beta+1)z\big)$.
Let $W(t)$ be a continuous lift such that $\exp(W(t))= \psi^{-1}(\gamma(t))$, defined for $t$ close enough to $1$. 
Then $f(\gamma(t))=h(W(t)) = \exp( (\beta+1) W(t))$. Since $\psi$ is an analytic diffeomorphism, $W(t)=L(t)+\on{const}+o(1)$ as $t\tend 1$.
Since $\big|\exp((\beta+1) W(t))\big| = \exp\big(\Re((\beta+1) W(t))\big) =\exp\big(\Re((1+i\tau) W(t))\big) =\exp\big(\Re W(t) - \tau\Im W(t)\big)$, the first claim follows.
\par
Second Claim: by homotopy, it is enough to prove it for paths contained in a neighborhood of $z_i$. In that case, we can apply the change of variable $\psi^{-1}$ as above, from which it follows that both continuations of $f$ tend to $0$.
\end{proof}

For instance, a curve tending to some $z_i$ and having a tangent there will satisfy the hypothesis of the lemma. 

\subsection{The parameter problem}

We will denote by $\log_p$ the principal branch of the logarithm: it is defined on $\C\setminus(-\infty,0]$ and takes values into the strip $\R\times(-\pi,\pi)$.

Recall that we have defined a multivaluated function
\[ m(z)=\left( \frac{(z-z_2)(z-z_4)}{(z-z_1)(z-z_3)} \right)^{\frac{\log K}{2i\pi}}.\]

\subsubsection{Theory}

If instead of the unit square $\sq$ on which the Beltrami form is constant, we decided to take a rectangle $R$ defined by $-h<\Re(z)<h$ and $-v<\Im(z)<v$, we could still define an associated affine surface and conformally uniformize it on the complement of four symmetric points in $\C$, with an isomorphism normalized at $\infty$ by $z\in \text{Chart }1\mapsto z+0+o(1)$. This would give another value of $z_1$, apart from which the affine charts would have exactly the same formula, with the same exponent $\log(K)/2i\pi$.
This allows to define a map
\[Z(h+iv) = z_1\]
that depends on $K$.

Independently of any choice of a rectangle, consider a given complex number $z_1$ in the upper right quadrant.
Consider the path $\gamma$ from $\infty$ to $z_1$ expressed by $z_1+(1+i)r$ with $r\in(0,+\infty]$ going from $+\infty$ to $0$.
Let
\[\Xi(z_1) = z + \int_{\gamma} (m(u)-1) du\]
where $m(z)$ denotes the continuous branch along $\gamma$ that equals $1$ at $\infty$. This map also depends on $K$.

If $z_1=Z(h+iv)$, then by Lemma~\ref{lem:sing} we must have $\Xi(z_1)=h+iv$: 
indeed let $\gamma_2$ be the image by $\phi$ of the path $h+iv+(1+i)r$, with $r\in(0,+\infty]$ going from $+\infty$ to $0$; first the path $\gamma$ is indeed homotopic in $\C\setminus\{z_1,z_2,z_3,z_4\}$ to $\gamma_2$, because both are contained in the upper right quadrant minus $\{z_1\}$; second the condition on the lifts $L(t)$ are satisfied: the lift of $\gamma$ is horizontal and the lift of $\gamma_2$ is asymptotic to a perpendicular to the line $\Re(z)-\tau\Im(z)=0$, as can be seen using Lemma~\ref{lem:loccorner}. In other words $\Xi\circ Z =\on{id}$ on the quadrant.
In particular, $Z$ is injective, $\Xi$ is surjective.
The map $Z$ and $\Xi$ are easily seen to be homogeneous: for all $z$ in the quadrant,
\[\forall \lambda>0,\quad Z(\lambda z)=\lambda Z(z) \quad\text{and}\quad \Xi(\lambda z)=\lambda \Xi(z).\]
Note that $Z(h+iv)$ is the image of the upper right corner of the rectangle $R$ by the straightening. It is thus a continuous function of $h+iv$, and it extends continuously on the boundary of the quadrant to the identity (because it then corresponds to a null Beltrami form). By this and homogeneity, $Z$ must be surjective, and tend to $\infty$ at $\infty$. The closed quadrant union $\infty$ is a compact set on which $Z$ is bijective, continuous and has an inverse $\Xi$. This implies that $\Xi$ is continuous. Therefore, $Z$ and $\Xi$ are inverse homeomorphisms.
Let us sum-up:
\begin{proposition}
The function $\Xi$ is a homeomorphism of the right quadrant, with a continuous extension to its boundary equal to the identity there.
\end{proposition}

\subsubsection{A method}

Here we present a numerical method to obtain $z_1$.
There are probably many other ways of doing this.

\begin{figure}[htbp]
  \begin{picture}(146.25,98.25)
    \put(0,0){\scalebox{0.75}{\includegraphics{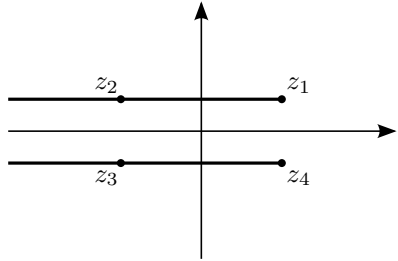}}}
    \put(105,65){$z_1$}
    \put(33,65){$z_2$}
    \put(33,30){$z_3$}
    \put(105,30){$z_4$}
  \end{picture}
  \caption{The domain of definition of $m_p$.}
  \label{fig:def}
\end{figure}

\begin{figure}[htbp]
  \begin{picture}(182.25,146.25)
    \put(0,0){\scalebox{0.75}{\includegraphics{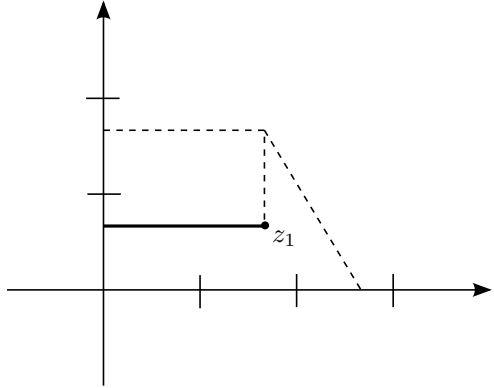}}}
    \put(100,55){$z_1$}
  \end{picture}
  \caption{The integration segments.}
  \label{fig:integralpath}
\end{figure}

To determine $z_1$, and thus the other $z_i$ by symmetry, we want to solve for \[ \Xi(z_1)=1+i.
\]
Let us choose a branch of the integrand function $m(z)$: let 
\[ m_p(z) = \exp\left(\frac{ \log K}{2i\pi} \Big( \log_p(z-z_2) + \log_p(z-z_4) - \log_p(z-z_1) - \log_p(z-z_3)\Big)\right)
\]
where $\log_p$ denotes the principal branch of the logarithm.
It is naturally defined on the complement in $\C$ of the two horizontal half-lines on the left starting from $z_1$ and $z_4$: 
\[ \on{Def}(m_p) = \C\setminus\big((z_1+(-\infty,0])\cup(z_4+(-\infty,0])\big).
\]
See figure~\ref{fig:def}.
Near $\infty$, it tends to $1$ and thus coincides with the branch used in Proposition~\ref{prop:expr1}, which is equal to $\phi^{-1}{}'$. Let us consider the continuation $f_p$ to $\on{Def}(m_p)$ of the affine chart $\phi^{-1}$ near $\infty$: $f_p(z)=z+\int_\infty^z (m_p(u)-1) du$ for any path of integration contained in $\on{Def}(m_p)\cup\{\infty\}$. Thus
\begin{equation}\label{eq:fpmp}
  f_p(z')-f_p(z)=\int_z^{z'} m_p(u) du
\end{equation}
for any path of integration from $z$ to $z'$ contained in $\on{Def}(m_p)$.
Note also that any path tending to $z_1$ within $\on{Def}(m_p)$ will satisfy the hypotheses of Lemma~\ref{lem:sing} and thus the formula defining $f_p(z)$ also works for $z=z_1$, and we have moreover
\[f_p(z_1)=\Xi(z_1).
\]
To better avoid numerical instabilities, especially when $\log K$ gets big, instead of computing $f_p(z_1)$ directly, in the author ended-up after a few attempts with the following way of determining $f_p(z_1)$. \emph{There is no claim that this choice is optimal or even good}, but it worked in the numerical applications.
Let $z_1=x+iy$, $z_2=-x+iy$, $z_3=-x-iy$ and $z_4=x-iy$. 
Compute
\begin{eqnarray*}
& & I_a=f_p(z_1+i)-f_p(iy+i), \\
& & I_b=f_p(z_1)-f_p(z_1+i), \\
& & I_c=f_p(z_1+i)-f_p(x+1), \\
\end{eqnarray*}
each of these numbers being determined by numerical integration of equation~\eqref{eq:fpmp} on the segment between the two points where $f_p$ is evaluated. See figure~\ref{fig:integralpath}.
Then
\[\Re(I_a+I_b)+i\Im(I_c+I_b)=\Xi(z_1).\]
So we have a method to compute $\Xi$. Recall that this function is homogeneous. A method to find the correct value of $z_1$, the one for which $\Xi(z_1)=1+i$, is then the following: 
\begin{itemize}
\item find the unique $\theta_1\in(0,\pi/2)$ such that the complex number $c=\Xi(\exp(i\theta_1))$ is collinear to $1+i$; for this we may look at the slope $s=\Re(c)/\Im(c)$ and use the secant method to solve equation $s=1$;
\item set $z_1 = \exp(i\theta_1)/\Re c$.
\end{itemize}

\subsubsection{Examples of parameter determinations}

To compute the integrals, the Simpson rule has been used, with an adaptive subinterval size, except near the singularity $z_1$, on which the integrand has been approximated by a constant times a complex power of the remaining distance, which can be integrated explicitly. Let us be more precise.

Recall that we have chosen a branch of the integrand function $m(z)$:
\[ m_p(z) = \exp\left(\frac{ \log K}{2i\pi} \Big( \log_p(z-z_2) + \log_p(z-z_4) - \log_p(z-z_1) - \log_p(z-z_3)\Big)\right)
\]
where $\log_p$ denotes the principal branch of the logarithm.

The adaptive method consists in subdividing an integration segment into unequal subintervals. The subintervals are constructed one after the other, in increasing order on the segment. The size of each subinterval is decided once its starting point is known. The error in Simpson's method, which consists in estimating $\int_a^b f$ by $\big(f(a)+4f((a+b)/2)+f(b)\big)/6$, is bounded from above by $(b-a)^5 \sup|f^{(4)}|/2880$. A very crude estimate was used for $\sup|m_p^{(4)}(z)|$ on $[a,b]$: $M\big|m_p(a)/|a-z_1|^4$ for some constant $M$. The integral on $[a,b]$ being close to $(b-a)|m_p(a)|$, the relative error is then estimated to $M\left|\frac{b-a}{a-z_1}\right|^4/2880$.
Then $b$ is chosen so that this quantity is small. This amounts to taking
\[b-a = d_0\times|a-z_1|\]
for some $d_0>0$.

For the segment that goes to $z_1$, we fix some distance $r_0>0$ to $z_1$ below which we stop using Simpson's rule, but instead do the following. Assume that we want to integrate $m_p(z)dz$ on the segment from $z_1+v$ to $v$. Depending on $v$, there is a constant $m_0\in\C^*$ such that $m_p(z_1+tv)\sim m_0 \exp(-\beta\log t))$ as $t\tend 0^+$, where $\beta=\log K/2i\pi$. The right hand side of this equivalence is easy to integrate on the segment from $z_1+v$ to $z_1$: $\int m_0 \exp(-\beta\log(t))) dz = vm_0\int_1^0 \exp(-\beta\log(t))dt = -v m_0/(1-\beta)$. In our case, $v=r_0 i$. Thus $m_0=\exp\Big(\beta\Big( \log_p(z_1-z_2) + \log_p(z_1-z_4) - \log_p(z_3-z_2) - \log_p(r_0 i) \Big)\Big)$.

In Table~\ref{tab:parampb} we show the results of some computer experiments. The values of $d_0$ and $r_0$ have been chosen by hand: the error has not been estimated, even if it is likely that the values are correct up to adding $\pm1$ to the last decimal. By convention all decimal numbers are approximations by truncation ($=$by default). The field ``steps in $\int$'' refers to the total number of subintervals the three integration intervals have been subdivided into. It does not count the subinterval that touches $z_1$. The values of $z_1$ and ``steps'' shown in the table are those corresponding to the last iteration of the secant method. The secant method never required more than 20 iterations, usually less (depending on the initial guess, obviously).

\begin{table}
\caption{}
\fbox{$\begin{array}{cccccc}
   K & \Re(z_1) & \Im(z_1) & d_0 & r_0 & \text{steps in }\int \\
   \hline
   2 & 1.2480 & 0.7676 & 0.2 & 0.01 & 38\\
  10 & 1.6483 & 0.4276 & 0.2 & 0.01 & 37\\
 100 & 1.8378 & 0.2350 & 0.2 & 0.005 & 39\\
 10^4 & 1.8976 & 0.1184 & 0.2 & 0.002 & 39\\
 10^6 & 1.9066 & 0.0789 & 0.2 & 0.001 & 46\\
 10^9 & 1.9104 & 0.0526 & 0.2 & 5.10^{-4} & 49\\
 10^{12} & 1.9117 & 0.0394 & 0.2 & 5.10^{-4} & 49\\
 10^{20} & 1.9127 & 0.0236 & 0.2 & 5.10^{-4} & 49\\
 10^{50} & 1.9132 & 0.00947 & 0.2 & 2.10^{-4} & 53\\
\end{array}$}
\label{tab:parampb}
\end{table}

As expected, $\Im(z_1)\tend 0$ as $K$ increases. It does so quite slowly: it seems to be comparable to $1/\log K$. It is not completely surprising as we know that in the formula for $m$ the constant $K$ appears only as $\log K$. As for $\Re(z_1)$ it seems to tend to a value close to $2$, but $<2$.

In these computations, the sufficient number of integration steps to reach a precision of $10^{-4}$ is surprisingly low. The same experiments have been run with lower values of $d_0$ and $r_0$, whence more integration steps, and yielded the same $10^{-4}$-approximations: on Table~\ref{tab:K2}, we see that the value stabilizes, and that $11$ decimal places could already be obtained with only approximately 1500 integration subdivisions.
\begin{table}
\caption{Numerical experiments for $K=2$.}
$\begin{array}{ccccc}
   \Re(z_1) & \Im(z_1) & d_0 & r_0 & \text{steps in }\int \\
   \hline
   1.25284529816 & 0.764992357071 & 1 & 0.1 & 5 \\
   1.33602299688 & 0.722976075222 & 1 & 0.001 & 5 \\   
   1.15536861038 & 0.460939047597 & 1 & 10^{-12} & 5 \\
   \hline
   1.23452767620 & 0.766441199359 & 0.5 & 0.5 & 9 \\
   1.24582003447 & 0.767669175930 & 0.5 & 0.2 & 10 \\
   1.24756645844 & 0.767688650200 & 0.5 & 0.1 & 11 \\
   1.24815951395 & 0.767643162006 & 0.5 & 0.01 & 14 \\
   1.24816586438 & 0.767640989544 & 0.5 & 10^{-4} & 21 \\
   1.24816586643 & 0.767640982471 & 0.5 & 10^{-8} & 34 \\
   \hline
   1.23447420199 & 0.766433209762 & 0.2 & 0.5 & 21 \\
   1.24748486405 & 0.767688303419 & 0.2 & 0.1 & 28 \\
   1.24807062425 & 0.767646071328 & 0.2 & 0.01 & 38 \\
   1.24807622667 & 0.767644491781 & 0.2 & 10^{-4} & 59 \\
   1.24807622705 & 0.767644491379 & 0.2 & 10^{-8} & 100 \\
   \hline
   1.24748382611 & 0.767688184430 & 0.01 & 0.1 & 564 \\
   1.24806951508 & 0.767645985170 & 0.01 & 0.01 & 793 \\
   1.24807511121 & 0.767644410897 & 0.01 & 10^{-4} & 1251 \\
   1.24807511157 & 0.767644410561 & 0.01 & 10^{-8} & 2167 \\
   \hline
   1.24748382610 & 0.767688184430 & 0.001 & 0.1 & 5633 \\
   1.24806951508 & 0.767645985169 & 0.001 & 0.01 & 7934 \\
   1.24807511120 & 0.767644410897 & 0.001 & 10^{-4} & 12537 \\
   1.24807511157 & 0.767644410560 & 0.001 & 10^{-8} & 21743 \\
   \hline
   1.24807511157 & 0.767644410560 & 10^{-5} & 10^{-8} & 2\,175\,079 \\
   \hline
   1.24807511157 & 0.767644410565 & 0.01 & 10^{-5} & 1480
\end{array}$
\label{tab:K2}
\end{table}
Similarly, for the extremely big value $K=10^{50}$, Table~\ref{tab:K1E50} shows that less than 2100 steps were enough to give an approximation with an error $<10^{-11}$.
\begin{table}
\caption{Numerical experiments for $K=10^{50}$.}
$\begin{array}{ccccc}
   \Re(z_1) & \Im(z_1) & d_0 & r_0 & \text{steps in }\int \\
   \hline
    1.91325679285 & 0.00947278823076 & 0.2 & 2.10^{-4} & 53 \\
    1.91325221530 & 0.00947281249592 & 0.2 & 10^{-8} & 97 \\
    1.91325221554 & 0.00947281249722 & 0.2 & 10^{-12} & 138 \\
    1.91325400121 & 0.00947285755832 & 0.1 & 10^{-8} & 202 \\
    1.91325400121 & 0.00947285755831 & 0.1 & 10^{-12} & 290 \\
    1.91325406829 & 0.00947285956006 & 0.05 & 10^{-8} & 414 \\
    1.91325407082 & 0.00947285966708 & 0.02 & 10^{-8} & 1045 \\
    1.91325407086 & 0.00947285966944 & 0.01 & 10^{-8} & 2098 \\
    1.91325407086 & 0.00947285966959 & 10^{-4} & 10^{-12} & 302810 
 \end{array}$
\label{tab:K1E50}
\end{table}
Again, in the absence of a abound on the error term, these experiments are given for information, not proof.

\subsection{Using the formula}

Recall that we have defined a multivaluated function
\[ m(z)=\left( \frac{(z-z_2)(z-z_4)}{(z-z_1)(z-z_3)} \right)^{\frac{\log K}{2i\pi}}.\]

\subsubsection{Computing the inverse of the straightening}\label{subsubsec:compphim}

Now, let us describe one way to determine whether a given point $z$ belongs to $\phi(\sq)$, and what is the value of $\phi^{-1}(z)$. Consider any path $\gamma$ from $\infty$ to $z$. Follow along this path the branch of $m(z)$ that equals $1$ at $\infty$. In order to follow $\phi^{-1}(z)$, begin by integrating $m(z)-1$ along the path. Then $\phi^{-1}(z)=z+\int_\infty^z (m-1)$ at the beginning, i.e.\ as long as $\phi^{-1}(\gamma(t))\notin \sq$. The fact that the path starts from $\infty$ is not a problem for the integration, because the integral is convergent near $\infty$ (if one prefers, applying the change of variable $z\mapsto 1/z$ leads to a standard path integral). Another possibility is to determine the value of $\phi^{-1}(z_0)$ for one point $z_0\in\C\setminus B(0,2)$, in which case $z_0\notin \phi(\overline{\sq})$, using either a path starting from $\infty$, or a power series expansion of the formula, and then use $z_0$ as a starting point for the path $\gamma$: $\phi^{-1}(\gamma(t))=\phi^{-1}(z_0)+\int_{z_0}^{\gamma(t)} m(z)dz$ as long as the computed value of $\phi^{-1}(\gamma(t))$ does not belong to $\sq$.
Whenever the computed value of $\phi^{-1}(\gamma(t))$ enters in $\sq$, the integration goes on but with an integrand equal to $A(m(z)dz)$ instead of $m(z)dz$, where $A(x+iy)=x+iKy$.
Whenever the computed value of $\phi^{-1}(\gamma(t))$ gets out of $\sq$, we switch back to the normal integrand. And so on, up to the end of $\gamma(t)$.
We could have started at $0$ instead of $\infty$, since we know that $\phi(0)=0$. In this case, begin by integrating $A(m(z)dz)$, and follow the branch of $m(z)$ that has the correct initial value (which is given in Proposition~\ref{prop:expr2}).

\subsubsection{Computing the the straightening}\label{subsubsec:compphid}

Computing $\phi(z)$ looks a lot like solving an Ordinary Differential Equation.
Recall that $\phi=\Phi\circ\Psi$.
Compared with the computation of $\phi^{-1}$, the difficult problem of determining whether a point $z\in\C$ belongs to $\phi(\sq)$ is replaced by the trivial problem of determining whether $z\in\sq$. So there is two cases: either $z\notin \sq$, in which case $\Psi(z)$ is the point with coordinates $w=z$ in Chart~1. Or $z\in\sq$ in which case $\Psi(z)$ is the point with coordinates $w=\Re(z)+i\Im(z)/K$ in Chart~2. 

Choose a point $w_0$ in the same chart for which $z_0=\Phi(w_0)$ and $\Phi'(w_0)$ are already known, choose any determination of $m(z_0)$ and choose a constant $a\in\C^*$ such that $1/\Phi'(w_0)=a\,m(z_0)$. For instance,\footnote{Another possibility is to let $w_0=\Phi^{-1}(z_0)$, for some chosen $z_0$, where $\Phi^{-1}$ and $m(z_0)$ may be computed with the method of section~\ref{subsubsec:compphim}.} in Chart~1 we may choose $w_0=\infty$, for which $z_0=\Phi(w_0)=\infty$ and $\Phi'(w_0)=1$, and choose $m(z_0)=1$ whence $a=1$. In Chart~2 we may choose $w_0=0$, for which $\Phi(w_0)=0$ and $\Phi'(w_0)=\exp(2\theta_1\log(K)/\pi)$ with the notations of Proposition~\ref{prop:expr2}, and choose $m(z_0) = \exp(-i 4 \theta_1 \log(K) / 2i\pi )$ whence $a=1$.

Choose a piecewise linear path $\gamma$ from $w_0$ to $w$ that does not leave the chart. Then
\[ \frac{\partial}{\partial t}\Phi(\gamma(t)) = \frac{\gamma'(t)}{ a\,m\big(\Phi(\gamma(t))\big)}
\]
where the correct determination of $m$ at $\Phi(\gamma(t))$ is obtained, as usual, by continuously following a branch along the part of the path $\Phi(\gamma(t))$ under construction.

\subsubsection{Drawing the boundary of the square}\label{subsubsec:db}

The previous method can be used to draw the boundary of the square: first determine $\phi(1)$ and $\phi(i)$, then draw the image of the path from $i$ to $1+i$, do the same for the path from $1$ to $1+i$, and use the symmetries to draw all the boundary of the square from these two curves.

However, this method is a bit unstable for drawing $\phi([i,1+i])$ when $K$ takes very big values. The reason is that the distance between $\phi(i)$ and $\phi(0)$ is comparable to $1/K$. (Also since $\phi(1)$ is only approximated, it means that the solution of the corresponding differential equation will miss $z_1$, even if slightly; since for big values of $K$ there is a lot of spiraling around $z_1$, could this yield big inaccuracies?)
The author ended up following the two paths in the opposite direction, i.e. from $1+i$ to $i$ and from $1+i$ to $1$. We cannot begin right on $1+i$, because $m$ degenerates at $z_1=\phi(i+1)$. So we need to begin slightly off: the difficulty is then to determine the initial values of $\Phi(z)$ and $\Phi'(z)$.

\begin{figure}%
\includegraphics[width=150pt]{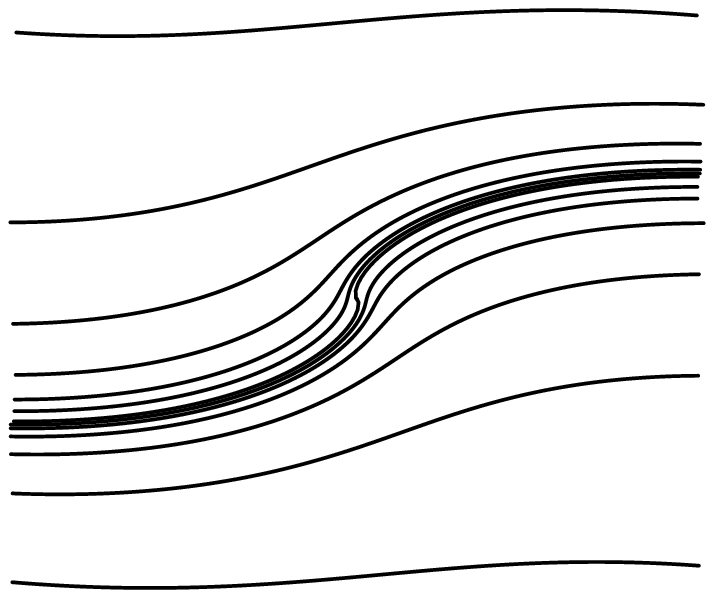}%
\includegraphics[width=150pt]{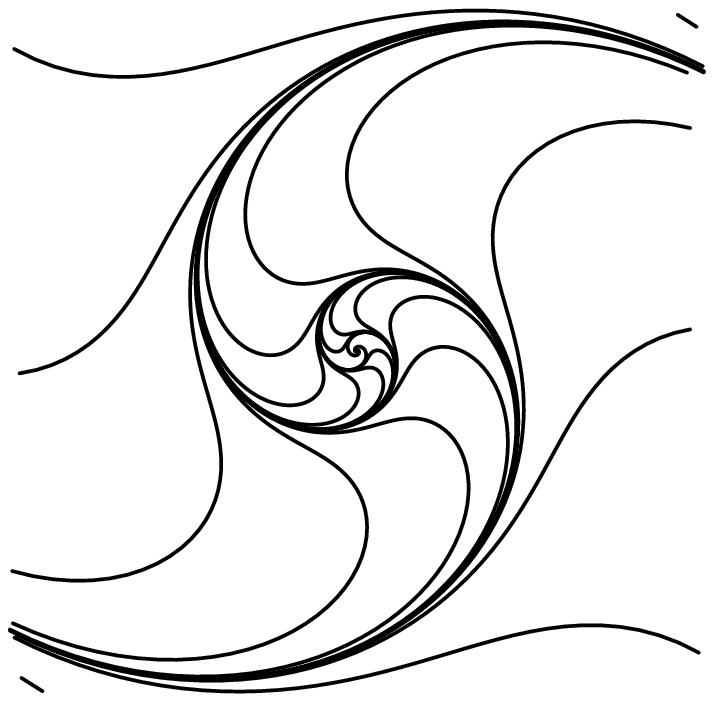}%
\caption{The field lines of $dz/dt=z^\beta$ for a purely imaginary $\beta$ are well defined as non-parameterized curves. They are shown for a small value of $\beta$ on the left, big on the right. The field lines passing close to the singularity stay close to the one passing through.}%
\label{fig:lf}%
\end{figure}

\begin{remark} It is not really important to start right on, or extremely close to, the spiraling image of the side of the square by $\phi$. Indeed, we are not at a fixed point of the vector field, but at another sort of singularity. Here small error is expected to remain small enough, see figure~\ref{fig:lf}. The problem comes from choosing the correct branch of $m(z)$. The fact that in our case $\beta$ is purely imaginary will make life easier.
\end{remark}

For this, let us use the fact that $\phi^{-1}$ is, along the curve $\phi([i,1+i))$, a primitive of a branch of $m(z)$. To make this branch more explicit, let $U=$ the upper right quadrant minus $\{z_1\}$. The restriction of $m$ to $U$ is multivaluated, and is well defined on the universal cover $\wt{U}$ of $U$. There, it satisfies
\[m(z)= \xi(z) \exp\big(\beta\log(z-z_1)\big)
\]
where $\beta = -\log(K)/2i\pi$ and where $\xi$ is a holomorphic function defined on the upper right quadrant, and 
\[\xi(z_1)=\exp\left( \frac{\log K}{2i\pi} \Big( \log_p(z_1-z_2) + \log_p(z_1-z_4) - \log_p(z_1-z_3)\Big)\right)
\]
where $\log_p$ denotes the principal branch of the logarithm.
A primitive $P$ of $m$ is a multivaluated function on $U$, well defined on $\wt{U}$, that satisfies
\[ P(z) =  \zeta(z) \exp\big((1+\beta)\log(z-z_1)\big) + \on{const}
\]
where $\zeta$ is a holomorphic function of $z$ defined on the upper right quadrant and
\[ c_0:=\zeta(z_1)=\xi(z_1)/(1+\beta) .
\]
Take $\on{const}=0$: then a point $z\in\wt{U}$ on the curve $\phi([i,i+1))$ satisfies $P(z)<0$. It turns out that, because $\beta$ is purely imaginary, the condition $P(z)<0$ is independent of the lift of $z\in U$ in $\wt{U}$, since changing the lift only multiplies $\exp\big((1+\beta)\log(z-z_1)\big)$ by $K^n>0$ for some $n\in\Z$.
From this it follows that the curve $\phi([i,i+1))$ spirals down to $z_i$ asymptotically on the logarithmic spiral defined by $z=z_1+r e^{i\theta}$ and
\[ \pi \equiv \arg(c_0) + \theta + \frac{\log K}{2\pi} \log_p (r) \pmod{2\pi}.
\]
One then starts from a point on this spiral for $r$ small enough and follows the differential equation 
\[ \frac{\partial z}{\partial t} = \frac{-1}{m(z)}
\]
for some branch of $m$, until the imaginary axis is crossed. 
\begin{remark} The time $t_1-t_0$ taken to reach this axis will not necessarily be $1$ but rather $K^n$ with $n\in\Z$, depending on which branch of $m$ is chosen. The author did not manage to find a simple way to determine, given the initial point $z_0=z_1+re^{i\theta}$, which branch of $m(z_0)$ to take in order to get $t_1-t_0=1$. However, it is likely that there is one.
\end{remark}

The same can be done for the curve $\phi([1,1+i))$, with
\[ -\frac{\pi}{2} \equiv \arg(c_0) + \theta + \frac{\log K}{2\pi} \log_p (r) \pmod{2\pi}.\]
and
\[\frac{\partial z}{\partial t} = \frac{-i}{m(z)}.\]

\subsubsection{Numerical experiments: a square on a diet}\label{sec:result}

\begin{figure}%
\begin{picture}(350,505)%
\put(175,385){%
\put(-55,10){\includegraphics[scale=0.5, bb=190 190 410 410]{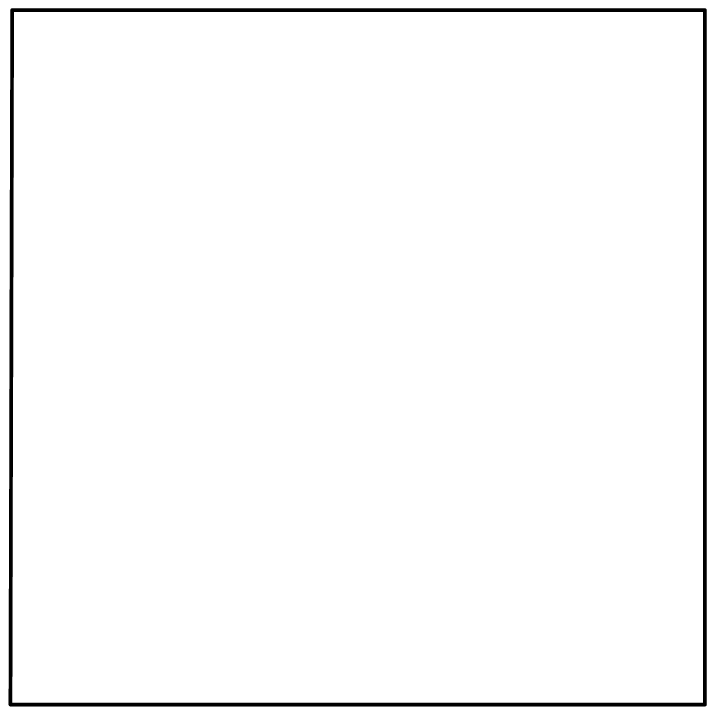}}%
\put(-13,63){$K=1$}
\put(-70,-80){\includegraphics[scale=0.5, bb=160 220 440 380]{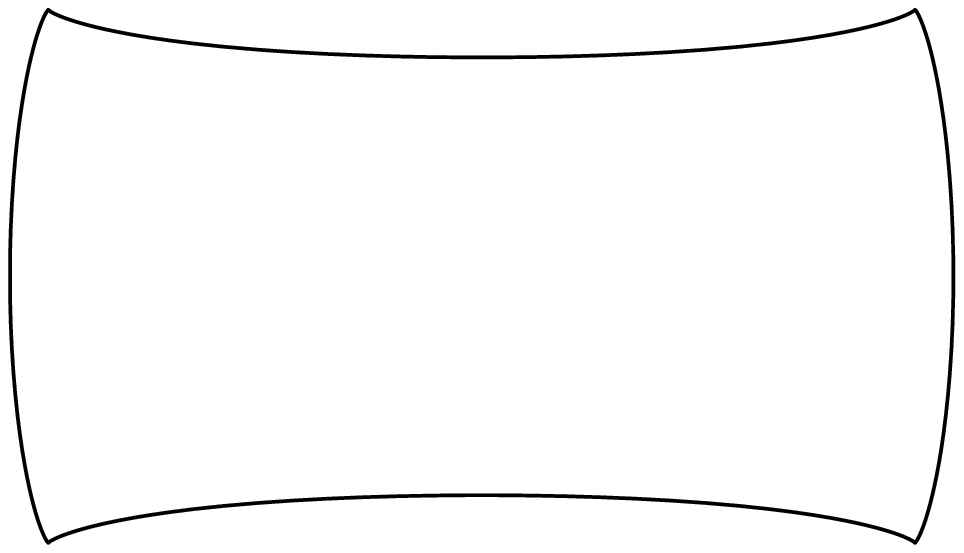}}%
\put(-0,-42){$2$}
\put(-87.5,-150){\includegraphics[scale=0.5, bb=125 240 475 360]{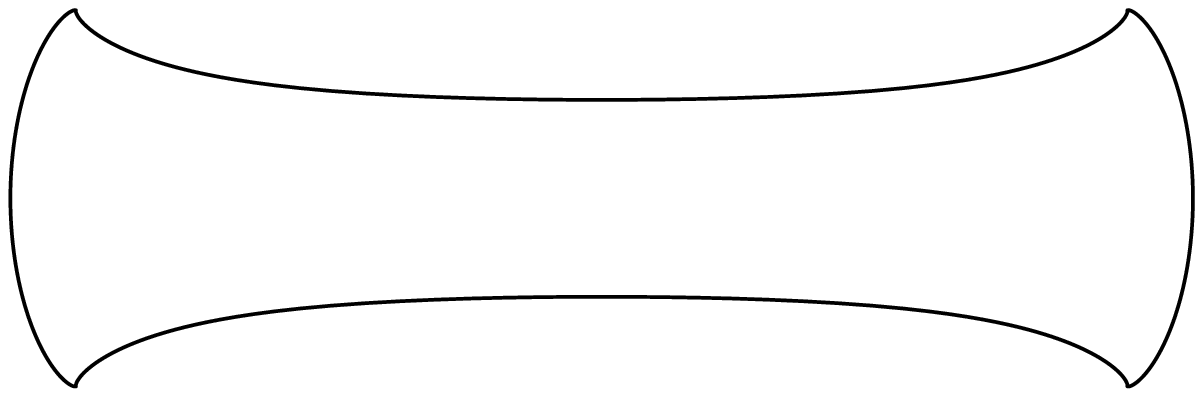}}%
\put(-0,-123){$5$}
\put(-100,-200){\includegraphics[scale=0.5, bb=100 260 500 340]{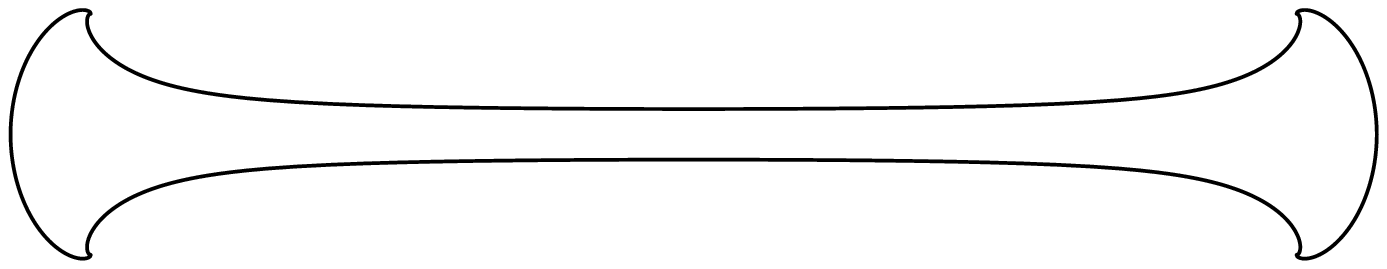}}%
\put(-0,-200){$20$}
\put(-110,-240){\includegraphics[scale=0.5, bb=80 270 520 330]{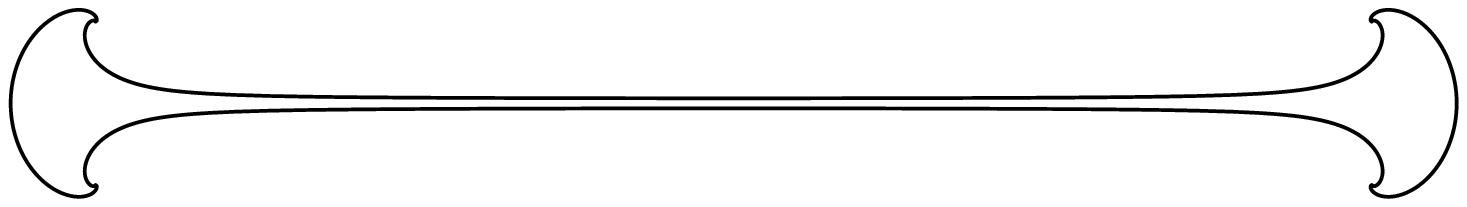}}%
\put(-0,-240){$100$}
\put(-110,-280){\includegraphics[scale=0.5, bb=80 270 520 330]{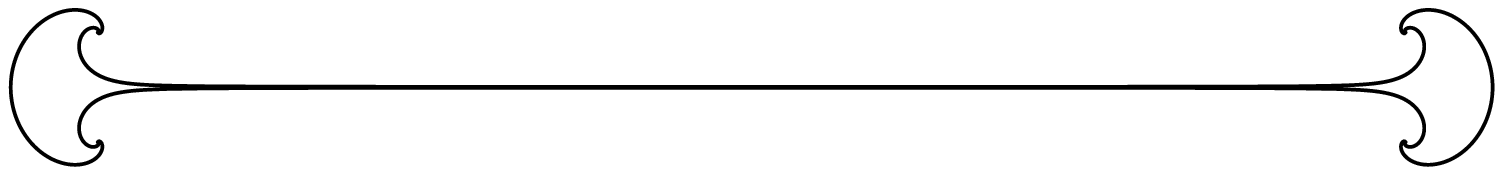}}%
\put(-0,-280){$1000$}
\put(-110,-315){\includegraphics[scale=0.5, bb=80 275 520 325]{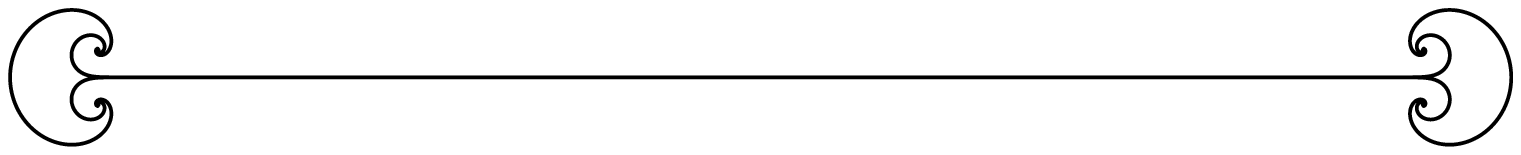}}%
\put(-0,-315){$10^6$}
\put(-110,-350){\includegraphics[scale=0.5, bb=80 275 520 325]{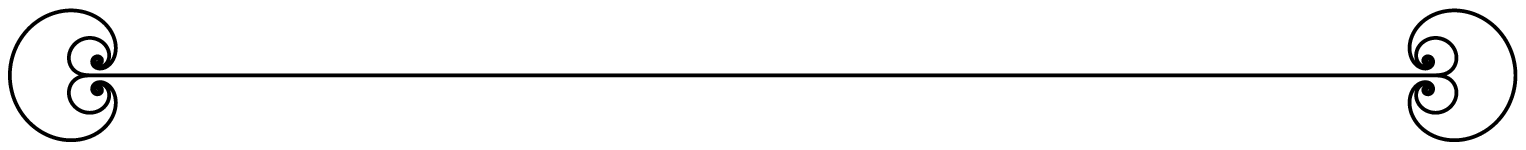}}%
\put(-0,-350){$10^{12}$}
\put(-110,-385){\includegraphics[scale=0.5, bb=80 275 520 325]{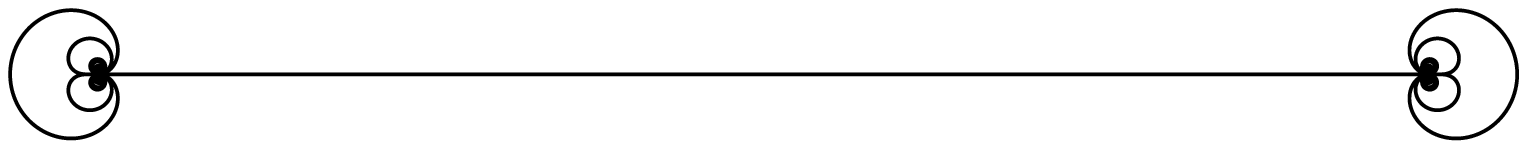}}%
\put(-0,-385){$10^{50}$}
}%
\end{picture}
\caption{The image of the square under the Beltrami form straightening, for different values of $K$.}%
\label{fig:sqseries}%
\end{figure}

\begin{figure}%
\setlength{\fboxsep}{0cm}%
\begin{picture}(300,450)%
\put(0,300){\includegraphics[scale=0.25]{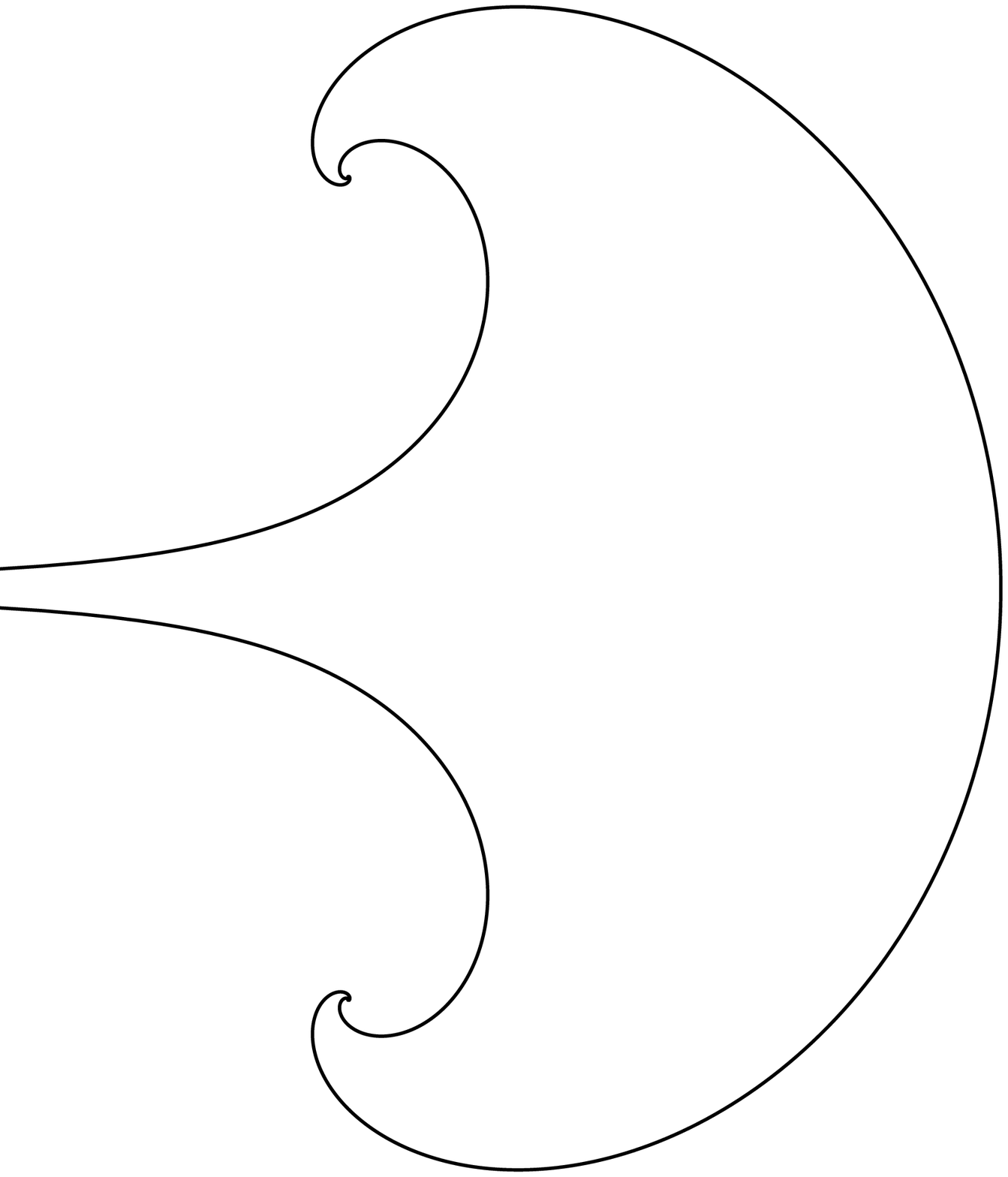}}%
\put(0,340){$K=1000$}
\put(150,300){\includegraphics[scale=0.25]{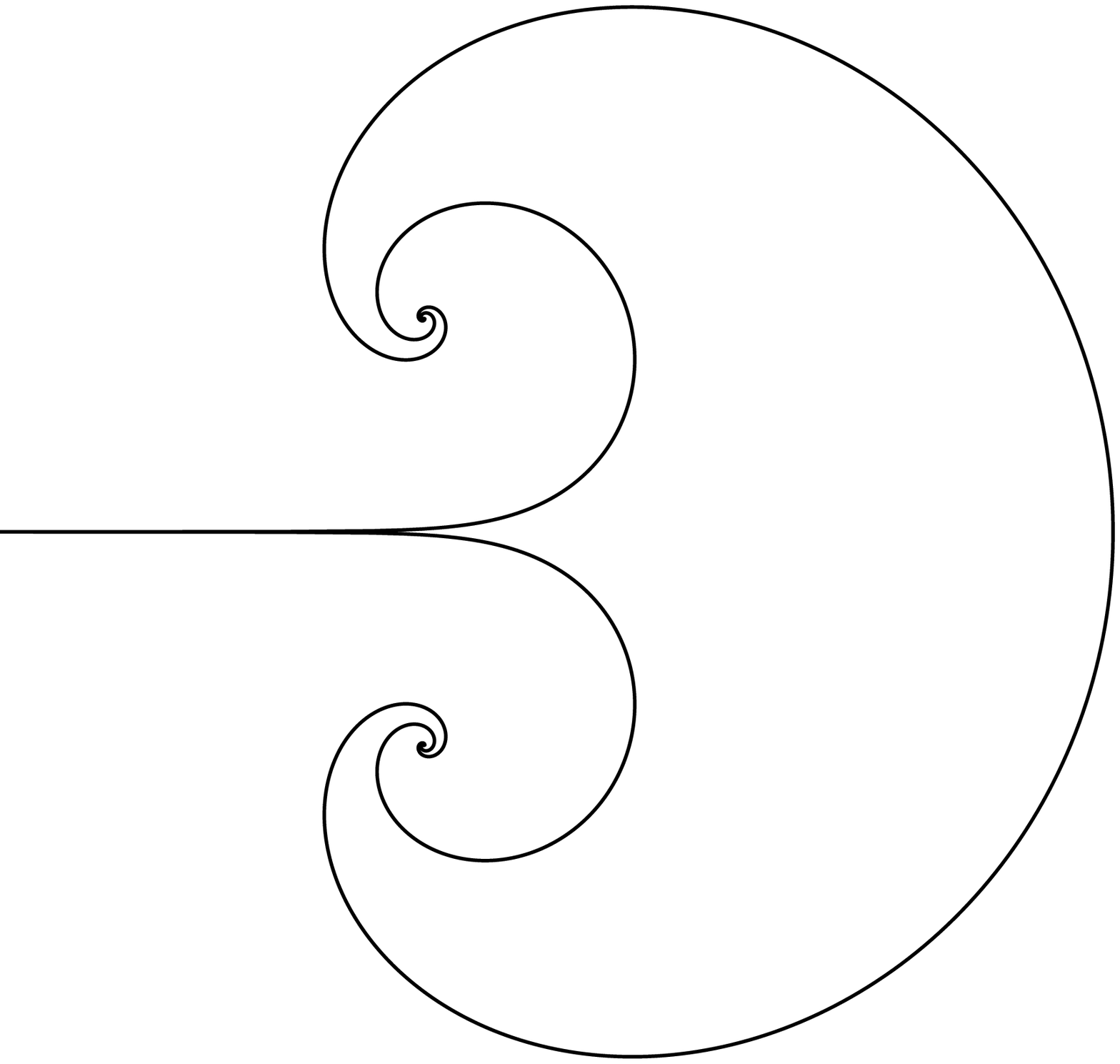}}%
\put(150,340){$10^6$}
\put(0,150){\includegraphics[scale=0.25]{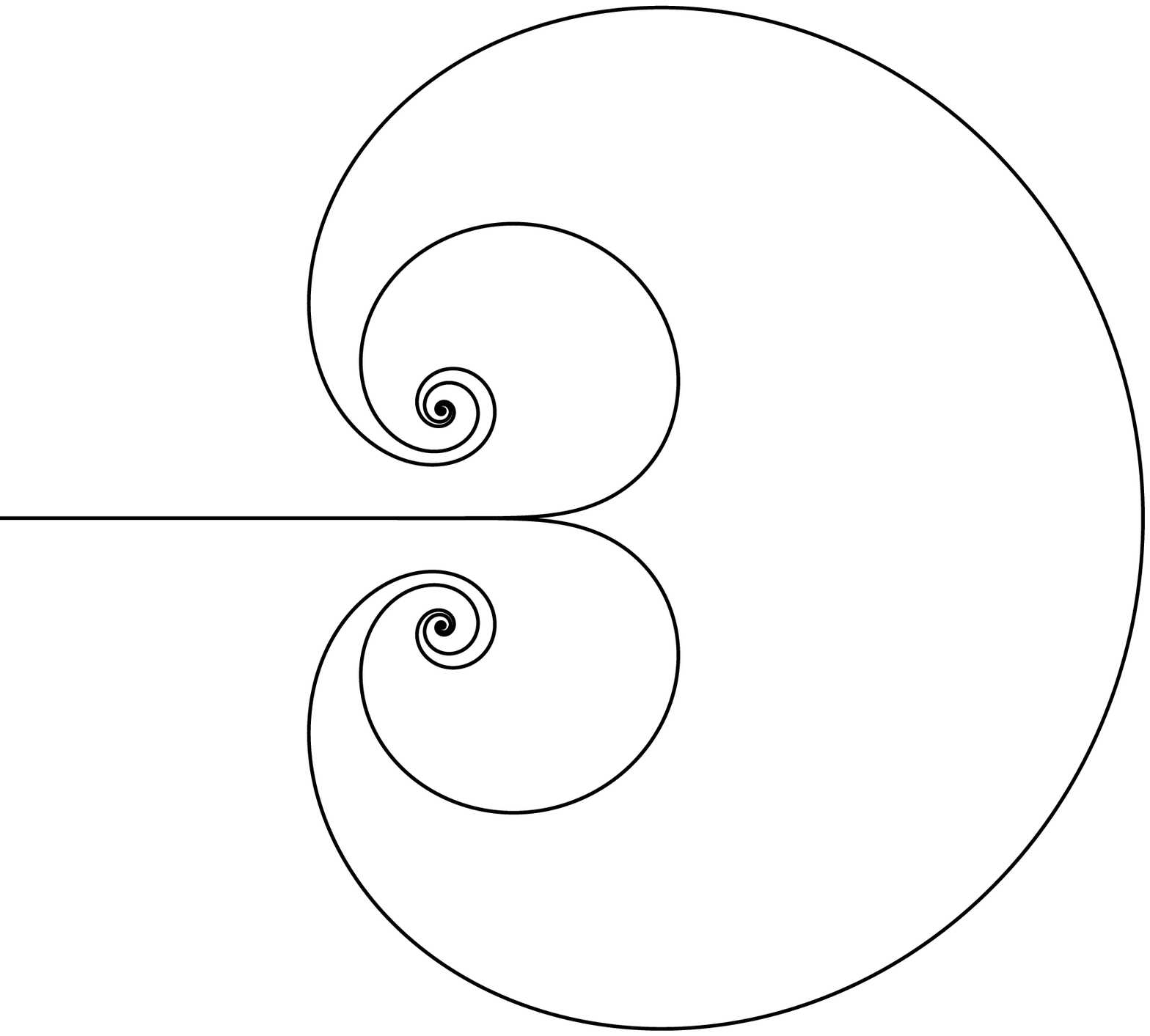}}%
\put(0,190){$10^{12}$}
\put(150,150){\includegraphics[scale=0.25]{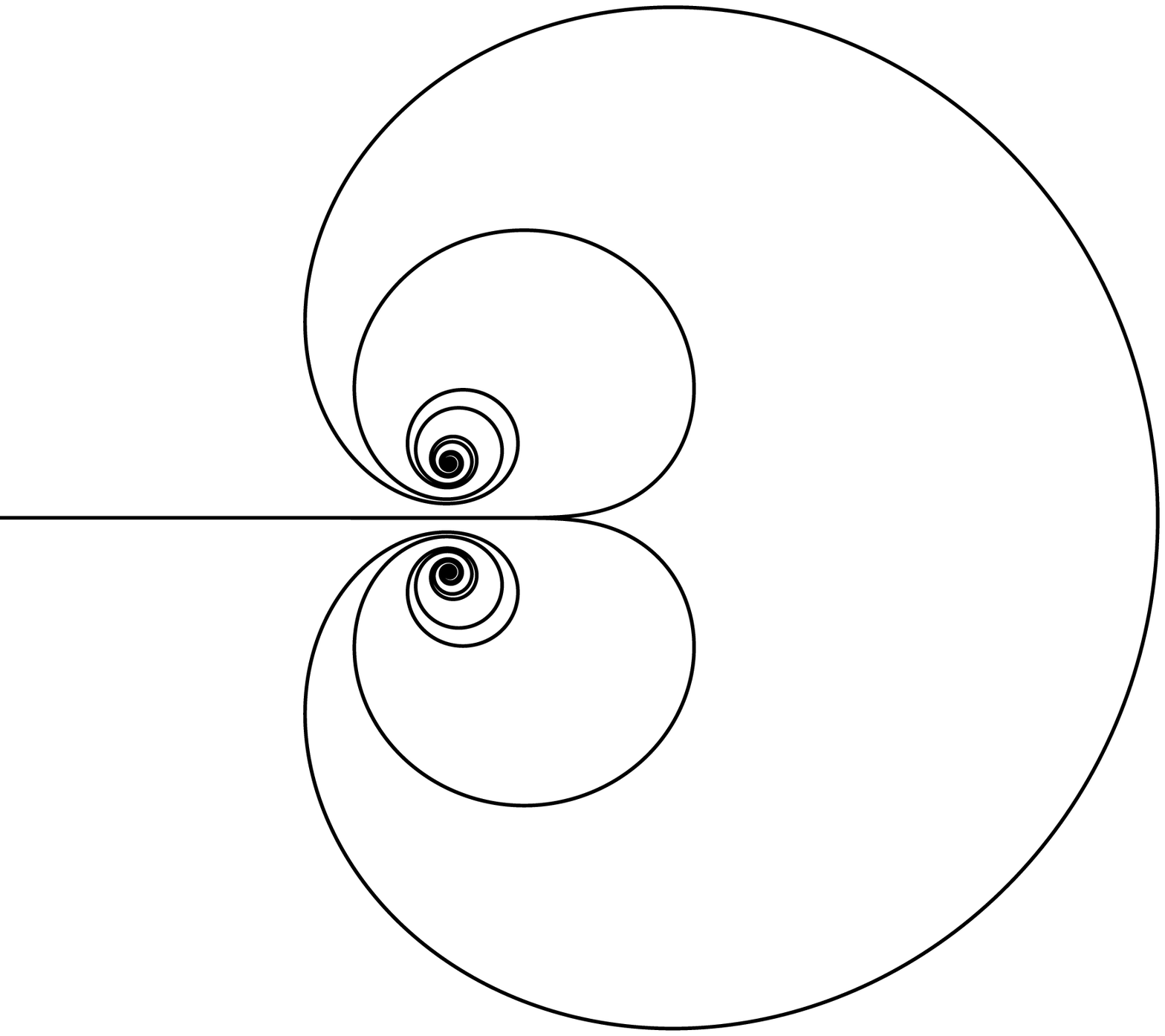}}%
\put(150,190){$10^{24}$}
\put(0,0){\includegraphics[scale=0.25]{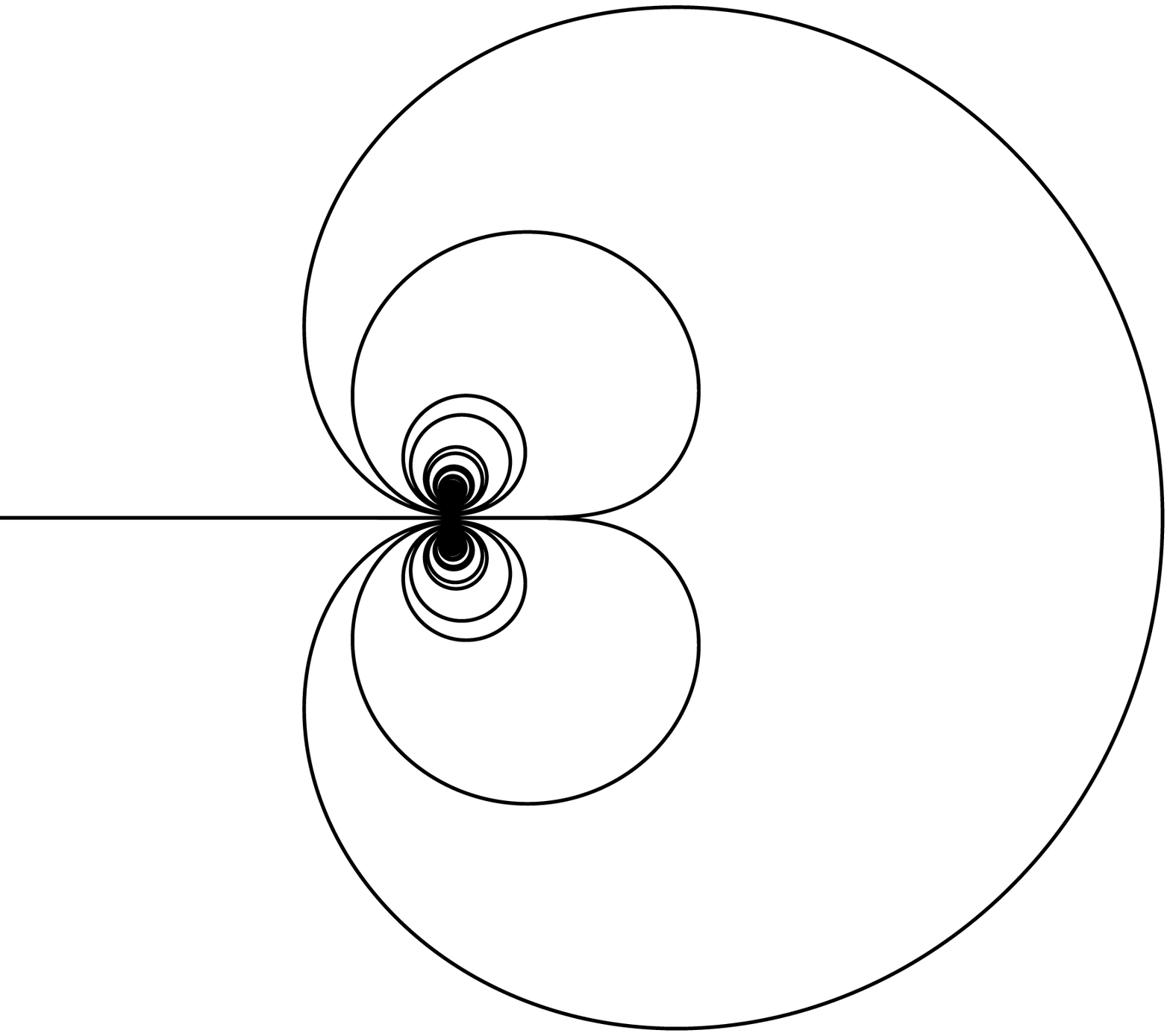}}%
\put(0,40){$10^{50}$}
\put(150,0){\fbox{\includegraphics[scale=0.25]{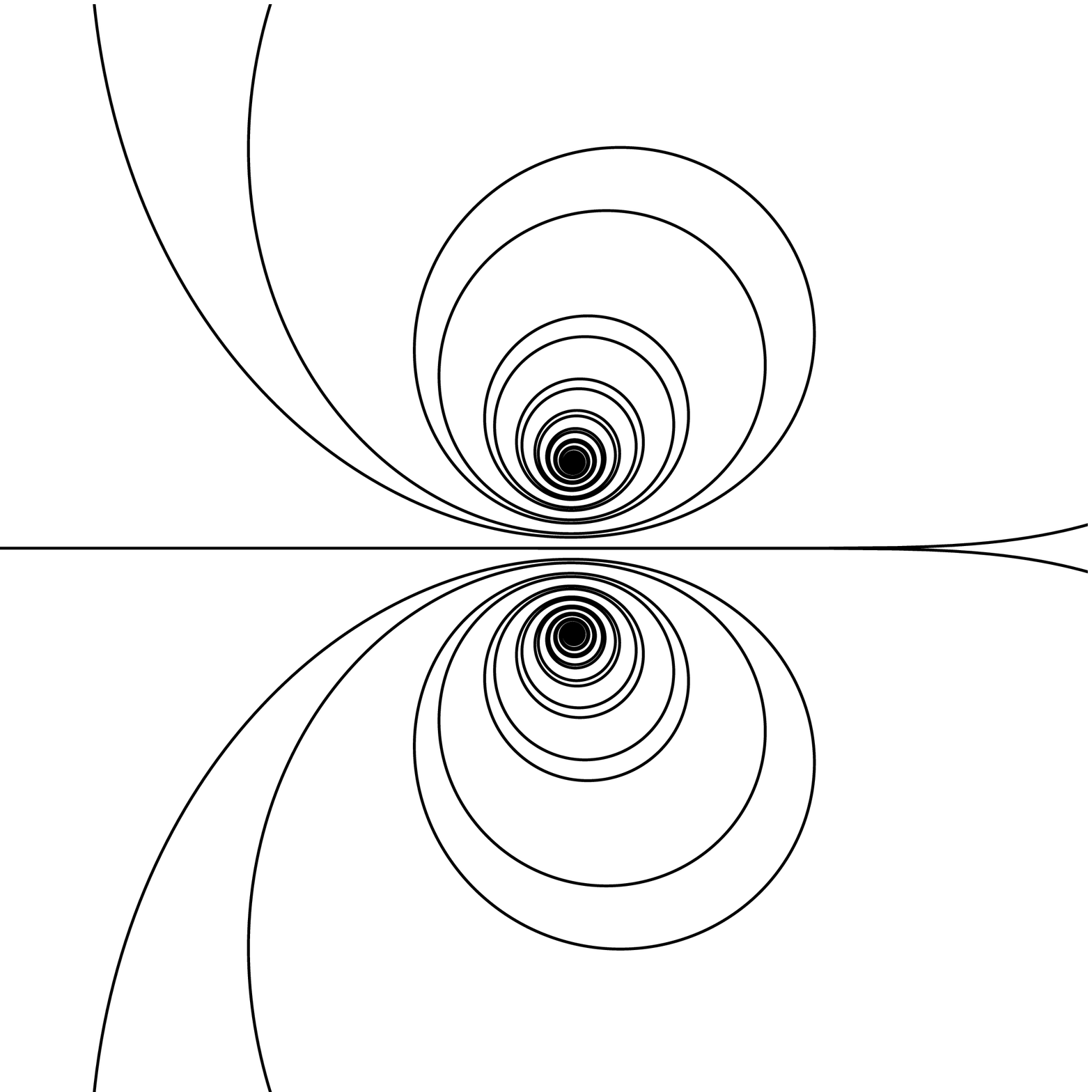}}}%
\put(200,5){$K=10^{50}$, closer view}
\end{picture}
\caption{Zoom on the bulb.}%
\label{fig:sqserieszoom}%
\end{figure}

We present on figures~\ref{fig:sqseries} and~\ref{fig:sqserieszoom} the image of the boundary of the square under the straightening $\phi$, computed using the method of section~\ref{subsubsec:db}. The differential equations $dz/dt=-1/m(z)$ and $dz/dt=-i/m(z)$ have been numerically solved with the fourth order Runge-Kutta method RK4, and an adaptive elementary interval size $dt$: more precisely $|dt|=\epsilon_0|m(z)|$, so that $|dz|$ is of order $\epsilon_0$. The value $\epsilon_0=10^{-5}$ was taken for the pictures. The list of points obtained approximates very well the curve, but is very big. A sublist was taken, and the corresponding polyline was drawn.

What we observe, as $K\tend+\infty$, is that:
\begin{itemize}
  \item $\Re(z_1) \tend x_\infty >0$
  \item $\Im(z_1) \tend 0$
	\item There seems to be a limit shape. 
	\item The limit shape is not a segment.
	\item It has non-empty interior.
	\item Yet it is not the closure of its interior: there is also a horizontal segment.
	\item The limit shape separates the plane into several connected components.
  \item The limit in the sense of Caratheodory of $\phi(\sq)$ (which is the union of the limit shape and all the bounded components of its complement), is equal to $[-x_\infty,x_\infty]\cup D \cup (-D)$ where $D$, is a topological disk (a Jordan curve and its bounded component) contained in ``$\Re(z)>0$'' and with $D\cap[0,x_\infty]=\{x_\infty\}$. It is smooth, except at $x_\infty$ where $D$ has a cusp. Let us call $D$ \emph{the bulb}.
	\item Near the endpoints of the segment, the boundary of the limit shape looks like a sheaf of tangent circles.
	\item Near these endpoints, when the parameter $K$ is big, the picture is reminiscent of parabolic bifurcation.
\end{itemize}

\section{The limit}\label{sec:limit}

\subsection{Behavior of the distortion derivative}

Let us note $z_1=x+iy$ with $x>0$ and $y>0$, so that
\bEA 
z_2=\phi(-1+i)=-x+iy, & & z_1=\phi(1+i)=x+iy, 
\\
z_3=\phi(-1-i)=-x-iy, & & z_4=\phi(1-i)=x-iy. 
\eEA
These numbers $x$ and $y$ depend on $K$. The previous sections suggest that, as $K\tend +\infty$, $y$ tends to $0$, $x$ tends to a positive value and $\phi$ has a non-affine limit on the complement of the square.

Whatever value of $K$ was chosen, the function $\phi$ is univalent on the complement of the square, thus in particular on $\C\setminus\overline{B}(0,\sqrt{2})$, and normalized so that $\phi(z)=z+0+o(1)$ as $z\tend\infty$, thus in particular $\phi'\tend 1$ at $\infty$. By a well known corollary of the area theorem $\phi(\C\setminus\overline{B}(0,\sqrt{2}))$ contains $\C\setminus\overline{B}(0,2\sqrt{2})$.
Therefore first
\[|z_i|\leq 2\sqrt{2};\] second $f=\phi^{-1}$ is univalent on $\C\setminus\overline{B}(0,2\sqrt{2})$ and normalized like $\phi$: $\phi^{-1}(z)=z+0+o(1)$ as $z\tend \infty$. For such normalized univalent functions $f$ there exists uniform bounds on the derivatives of a given order: in particular there exist functions $A_0(R)$, $A_1(R)$ and $A_2(R)$ independent of $f$ such that
\[|z|=R \implies A_0(R)\leq |f'(z)| \leq A_1(R)\text{ and } |f''(z)|\leq A_2(R).\]
In particular $|f''(z)/f'(z)|<A_2(R)/A_0(R)$.
Now, recall that
\[\frac{f''(z)}{f'(z)} = \eta(z) = \frac{\log K}{2i\pi}\left( \frac{1}{z-z_2} + \frac{1}{z-z_4} - \frac{1}{z-z_1} - \frac{1}{z-z_3}\right).\]
Thus for all fixed $z\in\C$ with $|z|>2\sqrt{2}$ one must have
\[\left( \frac{1}{z-z_2} + \frac{1}{z-z_4} - \frac{1}{z-z_1} - \frac{1}{z-z_3}\right)\underset{K\to +\infty}{\tend} 0.\]
Since the $z_i$ are bounded, this implies that collisions must occur: for all sequence $K_n\tend +\infty$ such that the four points $z_i$ converge, we have 
$\frac{1}{z-\lim z_2} + \frac{1}{z-\lim z_4} - \frac{1}{z-\lim z_1} - \frac{1}{z-\lim z_3} =0$ whenever $|z|>2\sqrt{2}$. Since the functions $1/(z-a)$ for different values of $a\in\C$ are linearly independent on any infinite subset of $\C$, then the following unordered pairs are equal: $\{\lim z_1,\lim z_3\} = \{\lim z_2,\lim z_4\}$.
Using the symmetries between the $z_i$, we conclude that
\[\text{either } x=\Re(z_1)\tend 0\text{ or } y=\Im(z_1)\tend 0\text{ or both}\]
as $K\tend +\infty$.
We can even be more precise: since
\[\left( \frac{1}{z-z_2} + \frac{1}{z-z_4} - \frac{1}{z-z_1} - \frac{1}{z-z_3}\right) = \frac{-4ixyz}{\prod_{i=1}^4(z-z_i)}\]
we have
\[xy \underset{K\to+\infty}{=} \cal O\Big(\frac{1}{\log K}\Big).\]
However, it seems hard to prove with this method what is suggested by the computer experiments of the previous sections: that $y$ tends to $0$ but not $x$, which shall tend to a positive value $x_\infty$, and that $\phi$ has a non-affine limit, thus $\eta$ has a non-zero limit, i.e.\ $xy\log K$ has a non zero limit.
This is the object of the next section.

\subsection{Behavior of the image of the corners when the dilatation ratio tends to $+\infty$.}

The objective of this section is summed-up in Theorem~\ref{thm:lims}.

\subsubsection{Motivating the estimation of some integrals}

Note that if the claims made above hold, then the function $\eta$ tends to
\[\eta(z)=\frac{\lim (y\log K)}{\pi}\left( \frac{-1}{(z-x_\infty)^2} + \frac{1}{(z+x_\infty)^2} \right).\]
This limit function defines an affine affine surface structure over $\C\setminus\{-x_\infty,x_\infty\}$, by expressing its distortion derivative. We may expect it to be a limit in some sense of the affine surfaces for $K<+\infty$. We may also look for a description of this new affine surface, in terms of charts and gluings.
From the pictures, it seems that the part corresponding to $\C\setminus D\cup(-D)$, where $D$ is the bulb, is obtained by taking the complement of the square and gluing the upper and lower sides by a translation.

Noteworthily, the function $m$ corresponding to the limit function $\eta$ would have no monodromy anymore:
\[m(z)=\exp\left(\frac{s}{(z-x_\infty)}-\frac{s}{(z+x_\infty)}\right)\]
with $s=\lim (y\log K)/\pi$.
A primitive of $m$ on the complement of $[-x_\infty,x_\infty]\cup D \cup (-D)$ would be a conformal map to the complement of the square.
Again from the pictures, and by analogy to the case $K<+\infty$ we expect that the pair $(x_\infty,s)$ is the unique pair of positive reals such that
\[\on{height}=\int_{t=+\infty}^0 (m(it)-1)d(it)=1\]
and
\[\on{width}=\int_0^1 m(t)dt=1.\]
The rigorous justification will be a bit complicated. We will begin with a slightly different version of the second integral, for $K<+\infty$, and let $K\tend+\infty$.

\subsubsection{Estimating integrals}

To estimate the limit of $z_1$ as $K\tend +\infty$, let us consider the continuation along the path from $+\infty$ to $0$ within $(0,+\infty)$ of the affine chart $\phi^{-1}$ near infinity. 
Because of the symmetries, the image by $\Phi^{-1}$ of the path ends up at the center of the square in Chart~2 (recall that $\Phi$ conformally maps the abstract affine surface to $\C\setminus\{z_1,z_2,z_3,z_4\}$).
The map 
\[z\mapsto \left\{ \begin{array}{ll}
\text{Chart}_1 \circ \Phi^{-1} & \text{on }\C\setminus\phi(\sq) \\
& \\
(1-K)+K\,\text{Chart}_2 \circ \Phi^{-1} & \text{on } \phi(\sq\cup (i-1,i+1))
\end{array}\right.\]
satisfies the equation $f''/f'=\eta$ and coincides with $\phi^{-1}$ out of the square. It is thus equal on $(0,+\infty)$ to 
the affine chart continuation mentioned above. Therefore
\[1-K = 0+\int_{+\infty}^0 (m(t)-1) dt\]
for the branch of $m$ on $(0,+\infty]$ that equals $1$ at $\infty$.
Let us make this more explicit:
\[1-K = \int_{t=+\infty}^0 \left( \exp\Big(\frac{\log K}{2i\pi} l(t)\Big)-1\right) dt\]
with
\[l(t)=\log_p(t-z_2) + \log_p(t-z_4) - \log_p(t-z_1) - \log_p(t-z_3)\]
where $\log_p$ denotes the principal branch of the logarithm. Using the symmetries between the $z_i$, this simplifies to
\[l(t)=i\Big(2\arg_p(t-z_2) - 2\arg_p(t-z_1)\Big)\]
where $\arg_p \in(-\pi,\pi)$.
And since $t>0$ and $\arg_p(t-z_2)\in (-\pi/2,0)$ and $\arg_p(t-z_1)\in(-\pi,0)$:
\[l(t)=i2\theta(t)\text{ with }\theta(t)=\arg_p\frac{t-z_2}{t-z_1} = \arg_p\frac{z_2-t}{z_1-t}.\]
Recall that
\[z_1=x+iy.\]
Then
\[\theta(t) = \arg_p \frac{-x-t+iy}{\phantom{-}x-t+iy} = -\arg_p (t^2+y^2-x^2 + i2xy) \in (0,\pi).\]
Thus
\bEA
\theta(t) & = & \cotan^{-1} X(t) \quad \text{with} \\
X(t) & = & \frac{t^2}{2xy} + \frac{1}{2}\Big(\frac{y}{x}-\frac{x}{y}\big)
\eEA
where $\cotan^{-1}$ takes its values in $(0,\pi)$.

So
\[K-1 = \int_{0}^{+\infty} \left( \exp\left(\frac{\log K}{\pi} \theta(t)\right)-1\right) dt = \underbrace{\int_0^{x} (\cdots) dt}_{\ds I_1}\ +\ \underbrace{\int_x^{+\infty} (\cdots) dt}_{\ds I_2}.\]

\titre{Bounding $I_2$} We will use the following inequality:
\[X>0 \implies \cotan^{-1}(X) < \frac{\pi/2}{1+X/\frac{\pi}{2}}.\]
For $I_2$, $t>x$ thus $X(t)\geq y/2x >0$, thus
\[I_2 < \int_{t=x}^{+\infty} \left( \exp\left(\frac{\log K}{2} \frac{1}{1+X(t)/\frac{\pi}{2}}\right)-1\right) dt. \]
Let us use the following convexity inequality:
\[0<s<S \implies e^s-1 <s\frac{e^S-1}{S}.\]
Applying this to $s=\frac{\log K}{2}/(1+X(t)/\frac{\pi}{2})$ and $S= \log(K)/2$ gives
\[I_2 < \int_{t=x}^{+\infty} \frac{\sqrt{K}-1}{1+X(t)/\frac{\pi}{2}} dt. \]
Thus
\bEA
  I_2 & < & \int_{t=x}^{+\infty} \frac{\sqrt{K}-1}{\ds 1+\frac{1}{\pi}\Big(\frac{t^2}{xy}+\frac{y}{x}-\frac{x}{y}\Big)} dt \ =\ \ldots
 \\
  \ldots & = & \int_{v=0}^{+\infty} \frac{\sqrt{K}-1}{\ds 1+\frac{1}{\pi}\Big(\frac{(v+x)^2}{xy}+\frac{y}{x}-\frac{x}{y}\Big)} dt
\ = \  \int_{v=0}^{+\infty} \frac{\sqrt{K}-1}{\ds 1+\frac{1}{\pi}\Big(\frac{v^2}{xy} +  \frac{2v}{y} +\frac{y}{x}\Big)} dt
  \\
  & < & \int_{v=0}^{+\infty} \frac{\sqrt{K}-1}{\ds 1+\frac{1}{\pi}\Big(\frac{v^2}{xy} +\frac{y}{x}\Big)} dt
        \ =\ \big(\sqrt{K}-1\big) \frac{\ds \pi\sqrt{\pi xy}}{\sqrt{1+\frac{y}{\pi x}}}
\eEA
because $\int_0^{+\infty} dx/(ax^2+b) = \frac{\pi}{2}/\sqrt{ab}$.
Now recall that $|z_i|\leq 2\sqrt{2}$, whence $x<2\sqrt{2}$ and $y< 2\sqrt{2}$.
Therefore
\[ I_2 < \big(\sqrt{K}-1\big)\pi\sqrt{8\pi}.
\]

\titre{Bounding $I_1$} It is a lot easier:
The integrand in $I_1$ is $<K-1$ and the integration interval is $(0,x)$ thus
\[I_1 < (K-1)x.\]
Now recall that $I_1+I_2=K-1$. Therefore
\[\liminf_{K\to+\infty} x \geq 1.\]
Let us use a slightly finer bound on the integrand: for $t\in(0,+\infty)$ we have $X(t) \in (X(0),+\infty)$ with $X(0) = \big(\frac{y}{x}-\frac{x}{y}\big)/2$.
Thus
\[ I_1 < (K^p-1)x < (K^p-1) 2\sqrt{2}\]
with $p=\frac{1}{\pi}\cotan^{-1} X(0) <1$. This implies that $p\tend 1$, i.e.\ that $\cotan^{-1} X(0) \tend \pi$, i.e.\ that $X(0) \tend -\infty$, i.e.\ that $y/x \tend 0$. Therefore
\[y \underset{K\to +\infty}{\tend} 0.\]

\titre{Making a pause} So far, we have proved that $0<\liminf x$ and $\lim y=0$. Recall that $|x|<2\sqrt{2}$ and that $xy\log K$ is bounded from above. We still need to prove that $x$ has a limit and that $xy\log K$ has a positive limit. Let us pass to a subsequence $K_n \tend+\infty$ such that $x$ has a limit $\chi$ and $xy\log K$ has a limit $\kappa$. Necessarily $\chi>0$. 

\begin{proposition}\label{prop:ints}
  \[\chi \int_0^1 \exp\left( \frac{-\sigma}{1-u^2} \right)du  = 1\]
with $\sigma= \frac{2\kappa}{\pi\chi^2}$, and
  \[\chi \int_0^{+\infty} \left(1-\exp\left( \frac{-\sigma}{1+u^2} \right) \right)du = 1\]
\end{proposition}

\titre{Proof of Proposition~\ref{prop:ints}}
\begin{lemma}\label{lem:I1}
  \[ \frac{I_1}{K_n} \underset{n\to\infty}\tend \chi \int_{0}^1 \exp\left(- \frac{2\kappa/\pi\chi^2}{1-u^2} \right) du
  \]
\end{lemma}
\begin{proof} We omit the index $n$.
  \[ \frac{I_1}{K} = \frac{-x}{K} + x\int_0^1 \exp\left(\frac{\log K}{\pi} \big(\theta(u)-\pi\big)\right)  du
  \]
  with
  \[\theta(u)=\cotan^{-1}\Big(\frac{u^2 x}{2y} + \frac{1}{2}\big(\frac{y}{x}-\frac{x}{y}\big)\Big).
  \]
  Now the integrand belongs to $(0,1)$ thus we may apply dominated convergence, and focus on the pointwise limit of the integrand for $n\tend+\infty$. Since $y\tend 0$ and $x$ is bounded from below, $y/x \tend 0$ and $x/y \tend +\infty$ thus for a fixed $u\in(0,1)$, $\frac{u^2 x}{2y} + \frac{1}{2}\big(\frac{y}{x}-\frac{x}{y}\big) \sim \frac{u^2-1}{2}\frac{x}{y} \tend -\infty$ whence  $\theta(u)=\pi + \frac{y}{x}\frac{2}{u^2-1} + o(y/x)$ whence $\frac{\log K}{\pi} \big(\theta(u)-\pi\big)\tend \frac{2\kappa/\pi\chi^2}{u^2-1}$.
\end{proof}
Since $K_n-1=I_1+I_2$ and $I_2=\cal O(\sqrt{K_n})$, Lemma~\ref{lem:I1} implies the first equality of Proposition~\ref{prop:ints}.

Let $Y$ denote the positive part of the imaginary axis: $Y=\setof{iy}{y>0}$.
Let $I_3=\int (m(z)-1) dz$ integrated along the path from $\infty$ to $0$ within $Y$, using the branch of $m$ defined on $Y$ and tending to $1$ at $\infty$. Then by similar arguments as for $I_1+I_2$, we have
\[I_3 = i - \frac{i}{K_n}\]
and
\[ I_3 = \int_0^{+\infty} \left(1- \exp\Big( \frac{\log K_n}{\pi} \on{arg}_p\big(x^2-y^2+t^2-i2xy\big) \Big)\right) i dt\]
where $\on{arg}_p$ denotes the principal argument $\in(-\pi,\pi)$.
\begin{lemma}
  \[ I_3 \underset{n\to\infty}\tend i\chi \int_{0}^{+\infty} \left(1-\exp\left(- \frac{2\kappa/\pi\chi^2}{1+u^2} \right) \right) du
  \]
\end{lemma} 
\begin{proof} Let $K=K_n$.
\[ I_3 = i\int_0^{+\infty} \left(1- \exp\Big( -\frac{\log K}{\pi} \cotan^{-1} \wt{X}(t) \Big)\right) dt\]
that we will write $I_3=i\int_0^{+\infty} f_n(t)dt$, with
\[\wt{X}(t) = \frac{t^2}{2xy} + \frac{1}{2}\Big(\frac{x}{y}-\frac{y}{x}\Big).\]
As $t$ varies from $0$ to $+\infty$, the function $\wt{X}(t)$ varies from $\frac{1}{2}\big(\frac{x}{y}-\frac{y}{x}\big)$ to $+\infty$. For $n$ big enough, $\frac{1}{2}\big(\frac{x}{y}-\frac{y}{x}\big)>0$. Then, using $1-\exp(-v)<v$ for $v>0$, and $\cotan v < 1/v$ for $v>0$, we get that for $n$ big enough, $\forall t>0$, $0<f(t)<\frac{2xy\log K}{(t^2+x^2-y^2)\pi}$. Since as $n\tend +\infty$, $xy\log K$ converges and $x^2-y^2 \tend \chi^2>0$, the sequence of integrals $I_3=i\int_0^{+\infty} f_n(t)dt$ is dominated. So we may focus on the pointwise limit of $f_n(t)$. The rest of the proof is similar to that of Lemma~\ref{lem:I1}, with the change of variable $t=xu$ and using $\cotan^{-1} v \sim 1/v$ as $v\tend +\infty$.
\end{proof}
This gives the proof of the second inequality in Proposition~\ref{prop:ints}.
\par\hfill\textit{End of the proof of Proposition~\ref{prop:ints}}\par\medskip

Let us recall the content of this proposition: for any sequence $K_n \tend +\infty$ such that $x$ has a limit $\chi$ and $xy\log K$ has a limit $\kappa$ we have
  \[\chi \int_0^1 \exp\left( \frac{-\sigma}{1-u^2} \right)du  = 1\]
with $\sigma= \frac{2\kappa}{\pi\chi^2}$, and
  \[\chi \int_0^{+\infty} \left(1-\exp\left( \frac{-\sigma}{1+u^2} \right) \right)du = 1\]
Now for a fixed $\chi$ the first integral is a decreasing function of $\sigma>0$ and the second one is increasing. By dominated convergence the first integral tends to $1$ as $\sigma\tend 0$ and to $0$ as $\sigma \tend +\infty$; the second tends to $0$ as $\sigma\tend 0$ and to $+\infty$ as $\sigma \tend +\infty$; both are continuous functions of $\sigma$. Therefore their quotient second$/$first, independent of $\chi$, is an increasing continuous bijection from $(0,+\infty)$ to itself. There is thus a unique value of $\sigma>0$ so that the quotient equals $1$. There is thus a unique pair $(\chi,\sigma)$ so that the pair of equations hold. There is thus a unique limit to $x$ and to $xy\log K$:
\begin{theorem}\label{thm:lims}
  1. There exists a unique solution $(x_\infty,\sigma)\in(0,+\infty)^2$ of the system
  \[ \left\{\begin{array}{rcl}
        W(x_\infty,\sigma) & = & 1
     \\ H(x_\infty,\sigma) & = & 1
     \end{array}\right.
  \]
  with
  \bEA
    W(x_\infty,\sigma) & = & x_\infty \int_0^1 \exp\left( \frac{-\sigma}{1-u^2} \right)du,
    \\
    H(x_\infty,\sigma) & = & x_\infty \int_0^{+\infty} \left(1-\exp\left( \frac{-\sigma}{1+u^2} \right) \right).
  \eEA
  2. As $K\tend +\infty$,
  \[z_1 \tend x_\infty, \quad \Im(z_1) \sim \frac{\pi\sigma x_\infty}{2\log K},\]
  \[\eta(z) \tend \frac{\sigma x_\infty}{2} \left( \frac{-1}{(z-x_\infty)^2} + \frac{1}{(z+x_\infty)^2} \right).\]
\end{theorem}

\begin{remark}
If instead of the square $\sq$ we had started from a rectangle of horizontal size $2w$ and vertical size $2h$, then the same conclusion holds with the system replaced by
  \[ \left\{\begin{array}{rcl}
        W(x_\infty,\sigma) & = & w
     \\ H(x_\infty,\sigma) & = & h
     \end{array}\right.
  \]
\end{remark}

\subsubsection{Numerical values}

Finding $\sigma$ and $x_\infty$ is similar to but simpler than finding $z_1$.
We obtain:

\bEA
 \lim \Re(z_1) = x_\infty & = & 1.9133480795\ldots
 \\
   \sigma & = & 0.3628700796\ldots
 \\
  \frac{\sigma x_\infty}{2} & = & 0.3471483850\ldots
 \\
 \lim \big(\Im(z_1)\log K\big) = \frac{\pi\sigma x_\infty}{2} & = & 1.0905988161\ldots
\eEA

For instance, if we take $K=10^{50}$, and the value computed for $z_1$ in table~\ref{tab:K1E50} then we obtain $\Re(z_1)=1.9132540708\ldots$ and $\Im(z_1)\log K = 1.0906032731\ldots$ which are close to the values above.

\subsubsection{A remark about a residue}

The following integral appearing in the second equation of theorem~\ref{thm:lims} can be expressed as a residue:
\[H=H(x_\infty,\sigma)=x_\infty \int_0^{+\infty} \left(1-\exp\left( \frac{-\sigma}{1+u^2} \right) \right)du\]
Recall that at the limit the function $m$ is not anymore multivaluated:
\[m(z)=\exp\left(\frac{s}{(z-x_\infty)}-\frac{s}{(z+x_\infty)}\right)\] with $s=\lim(y\log K)/\pi=\sigma x_0 /2$.
Recall that $iH=\int(m(z)-1)dz$, integrated along the path from $+\infty$ to 0 in the imaginary axis. By symmetry $2iH=\int(m(z)-1)dz$ along the full imaginary axis followed downwards.
Consider the integral $\int(m(z)-1)dz$ along the following curve: it starts from $iR$ and descends the imaginary axis down to $-iR$; then it follows the half-circle from $-iR$ back to $iR$ and located on the right of the imaginary axis. This integral will be equal to $2i\pi$ times the residue of $m$ at $x_\infty$ on one hand. On the other hand, the vertical part tends to $2iH$ as $R\tend +\infty$, and the circular part is $\cal O(1/R)$ because $m(z)=1 + \cal O(z^2)$ as $z\tend\infty$. Therefore
\[H=\pi \times\text{residue(}m\text{ at }x_\infty).\]
However, this residue is not really explicit. It can be computed using power series expansion and gives
\[\on{residue} = \frac{x_\infty \sigma}{2} \sum_{n=0}^{+\infty} \left(\frac{-\sigma}{4}\right)^n \frac{(2n)!}{(n!)^2(n+1)!}\]
which converges quite rapidly (like the power series of $\exp(z)$).

But in fact, the same expression could have been computed directly from the definition of $H$ as an integral: by expanding $1-\exp\Big(\frac{-\sigma}{1+u^2}\Big) = -\sum_{n\geq 1} \frac{1}{n!} \Big(\frac{-\sigma}{1+u^2}\Big)^n$, integrating from $0$ to $\infty$, permuting the sum and the integral (it is possible because the total mass is bounded), and shifting the index, we get
\[H=x_\infty \sigma\sum_{n=0}^{+\infty} \frac{(-\sigma)^n}{(n+1)!}\int_0^{+\infty} \frac{du}{(1+u^2)^n}\]
and the constant $\int_0^{+\infty} \frac{du}{(1+u^2)^n}$ can be computed explicitly, by induction (or as a residue\ldots).

The author ignores if an analog approach can be set up for the integral $W(x_\infty,\sigma)$.

\subsection{The limit affine surface}

\subsubsection{Introduction} Recall that any analytic function $h$ on an open subset $U$ of $\C$ defines an affine surface structure on $U$, whose charts are the local solutions of $f''/f'=h$ (in particular they are analytic). We called $h$ the distortion derivative of the affine structure.

The function $\eta$, which depends on $K$, is defined as the distortion derivative of some affine structure on $\C\setminus\{z_1,z_2,z_3,z_4\}$, which was obtained by conformally uniformizing an affine surface $\cal S$ defined by polygonal pieces in $\C$ glued along their edges by affine maps.
It is thus natural to interpret its limit as $K\tend +\infty$:
\[\eta_\infty(z) = \frac{\sigma x_\infty}{2}\left( \frac{-1}{(z-z_\infty)^2}+\frac{1}{(z+z_\infty)^2}\right)\]
as the distortion derivative of an affine structure on $\C\setminus\{-z_\infty,z_\infty\}$.

It would be nice to:
\begin{itemize}
\item have an isomorphic model of this affine surface as polygonal pieces in $\C$ glued along their edges by affine maps,
\item precise in which sense it can be considered as a limit of the sequence of affine surfaces $\cal S$, independently of their conformal uniformization to $\C\setminus\{z_1,z_2,z_3,z_4\}$,
\item reprove the convergence of the corresponding function $\eta$, independently of the previous sections, i.e.\ without computing and estimating all these integrals; in particular, the uniqueness of the extracted limits of $\Re(z_1)$ and $\Im(z_1)\log K$ would follow from the uniqueness of (normalized) conformal uniformization.
\end{itemize}
In this section, we will present a candidate isomorphic model for the limit affine surface.

\subsubsection{The candidate}

Let us motivate the definition of this candidate.

\begin{figure}[htbp]%
\begin{picture}(350,400)
\put(0,0){\includegraphics[width=350pt]{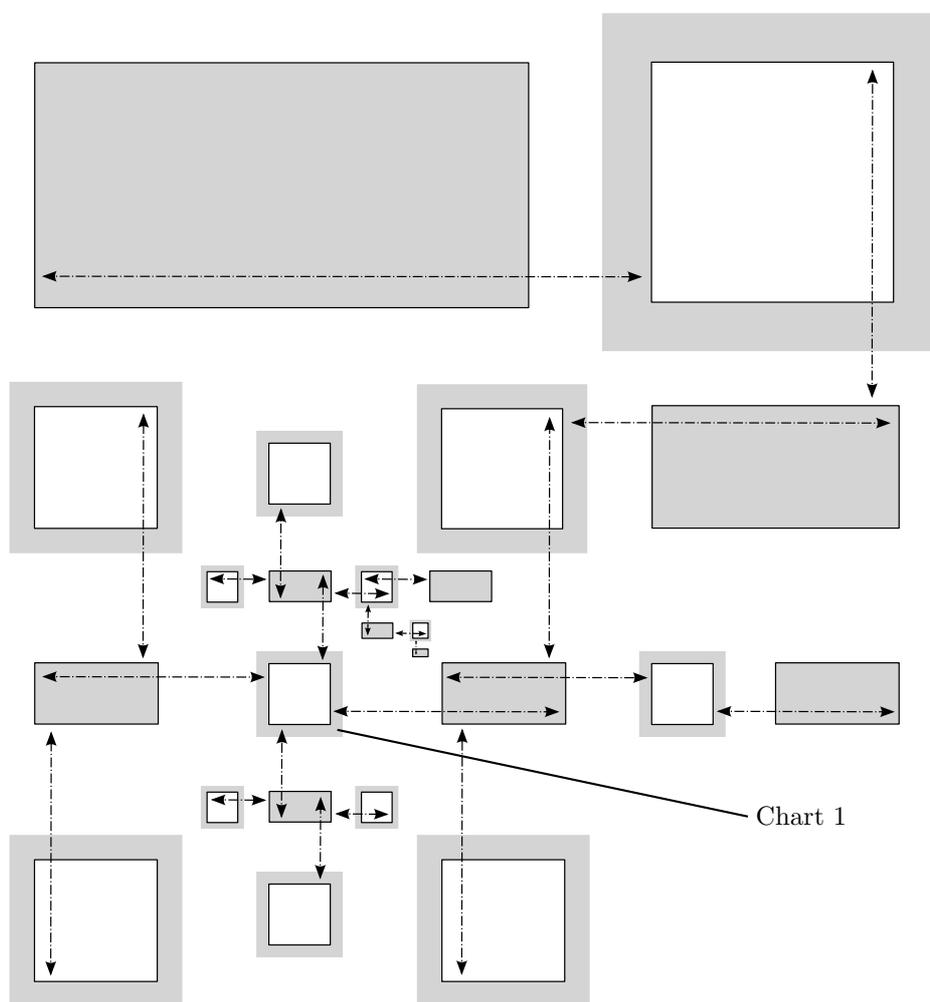}}%
\put(280,69){Chart 1}
\put(277,72){\psline(-5.4,1.15)}
\end{picture}
\caption{Part of the construction of the Riemann surface over $\C$ corresponding to the developing map of $\cal S$ with Chart~1 as a germ.  The picture illustrates the case $K=2$. The arrows indicate which sides shall be glued together.}%
\label{fig:ucissm}%
\end{figure}

Given an affine surface $\cal S$, the affine charts can be lifted to its universal cover $\wt{\cal S}$ and define an affine surface structure on it. Moreover, since $\wt{\cal S}$ is simply connected, any germ of an affine chart extends uniquely to all of $\wt{\cal S}$. This extension is called the \emph{developing map}. It is not necessarily injective anymore. This developing map $\wt{\cal S}\to\C$ can be seen as a structure of ``Riemann surface over $\C$''. Do not be mislead by the terminology: a Riemann surface over $\C$ is much more rigid than a Riemann surface; the transition maps between charts are required to be the identity; so it is more rigid than an affine surface or a translation surface.

For a given $K<+\infty$, let us look at a little variant: consider the usual affine surface $\cal S$, which we recall is defined affine gluings of two pieces: the complement of a square (Chart~1) and a rectangle (Chart~2). Instead of looking at the universal cover of $\cal S$, let us look at the universal cover of $\cal S \cup\{\infty\}$. Remove to this universal cover all points that project to $\infty$. Call this affine surface $\cal S'$. It has also a well defined developing map, because such a map can be define by gluing together, by $z\mapsto z$, affine images of the square's complement and of the rectangle. This is illustrated for $K=2$ in figure~\ref{fig:ucissm}.

\begin{figure}[htbp]%
\begin{picture}(350,290)
\put(0,0){\includegraphics[width=350pt]{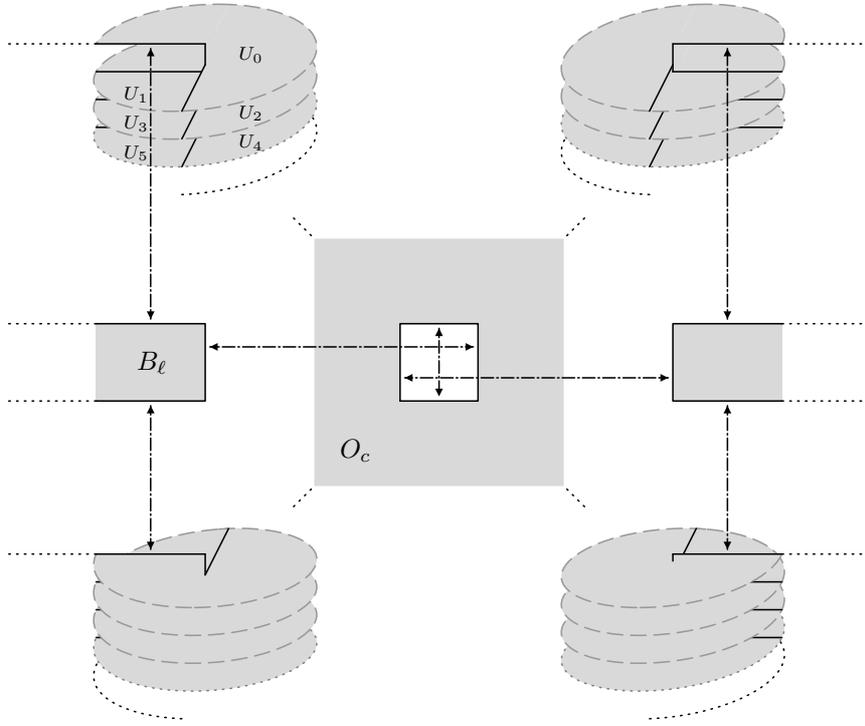}}
\put(138,110){$O_c$}
\put(62,143){$B_\ell$}
\put(100,260){$\scriptstyle U_0$}
\put(57,245){$\scriptstyle U_1$}
\put(100,238){$\scriptstyle U_2$}
\put(57,234){$\scriptstyle U_3$}
\put(100,227){$\scriptstyle U_4$}
\put(57,223){$\scriptstyle U_5$}
\end{picture}%
\caption{A flat model of the limit affine surface. The arrows indicate which side shall be glued together. The spiraling stacks represent infinitely many copies of the slit plane $\C\setminus(-\infty,0]$, glued together along the slit by $z\mapsto z$; their affine structure is the one induced by the inclusion in $\C$.}%
\label{fig:RS}%
\end{figure}

\begin{figure}[htbp]%
\begin{picture}(350,234)
\put(0,0){\includegraphics[width=350pt, bb=0 100 600 500]{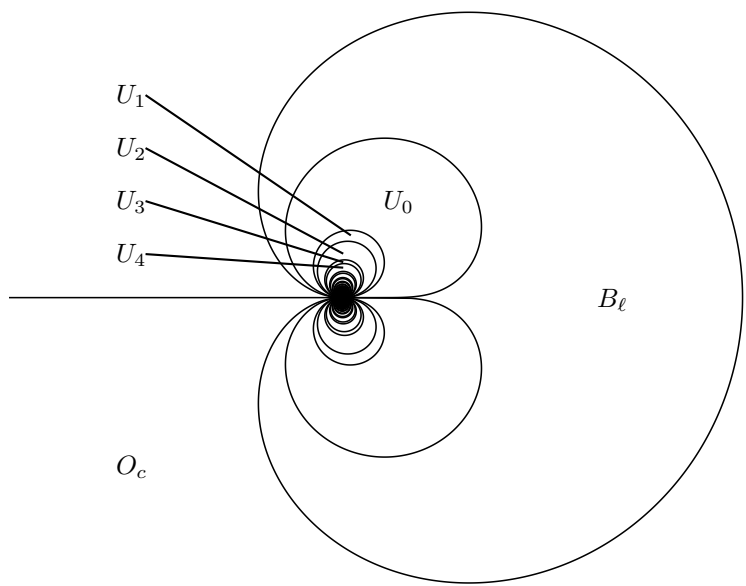}}
\put(40,50){$O_c$}
\put(220,113){$B_\ell$}
\put(140,150){$U_0$}
\put(40,190){$U_1$}
\put(51,193){\psline(2.7,-1.86)}
\put(40,170){$U_2$}
\put(51,173){\psline(2.6,-1.4)}
\put(40,150){$U_3$}
\put(51,153){\psline(2.6,-0.81)}
\put(40,130){$U_4$}
\put(51,133){\psline(2.6,-0.18)}
\end{picture}%
\caption{The corresponding regions, after uniformization.}%
\label{fig:lim}%
\end{figure}

Playing with this, letting $K\tend +\infty$, we end up suspecting that the following surface is a natural limit of the affine surface, adapted to the normalization at $\infty$ (a.k.a.\ $\Phi'(z)=z+0+o(1)$ in Chart~1).
It is defined as follows, and is illustrated in figure~\ref{fig:RS}. Let $O_c=\C\setminus\sq$. Let $B_\ell$ be defined by ``$\Re z < 1 \ \&\ |\Im z|<1$''. It is seen as a limit of the rectangle, whose width went to $+\infty$ on the left. Similarly, let $B_r$ be defined by ``$\Re z> -1 \ \&\ |\Im z|< 1$''. Glue $B_\ell$ to $O_c$ along the open segment $(1-i,1+i)$ (recall this means that we thicken a little bit the open sets along this segment and glue them by $z\mapsto z$). Glue $B_r$ to $O_c$ along $(-1-i,-1+i)$.  Glue the upper side and the lower side of the square bounding $O_c$, ends excluded, by a translation.  Now we will add four chains of pieces on respectively the upper side of $B_\ell$, its lower side, the upper side of $B_r$ and its lower side. We will describe only the first chain, the others are symmetric. Let $U_{2n} = \setof{z\in\C}{\Re z>0\text{ or }\Im z>0}$ and $U_{2n+1}=\setof{z\in\C}{\Re z<0\text{ and }\Im z<0}$. Glue $U_0$ to $B_\ell$ along $(-\infty,0)$, $U_{2n}$ to $U_{2n+1}$ along $(-\infty,0)\times i$ and $U_{2n+1}$ to $U_{2n+2}$ along $(-\infty,0)$. Gluing together the $U_n$ gives a piece that can also be seen as the subset ``$\arg(z)\in(-\infty,\pi)$'' of the universal cover of $\C^*$ on which we would have put the following affine structure: the natural affine structure on $\C$ is first pulled back on $\C^*$ by the inclusion map, then pulled-back by the universal cover of $\C^*$. This is represented as spirals in figure~\ref{fig:RS}. The uniformization to $\C\setminus\{x_\infty,-x_\infty\}$ is supposed to yield figure~\ref{fig:lim}, which was drawn using the formula for the limit distortion derivative $\eta_\infty$ and drawing the curves corresponding to the black lines in figure~\ref{fig:RS}.

\section{Conclusion}\label{sec:ccl}

The story is not finished. We may wonder what is the limit when the major axis of the ellipses are slanted. The same approach as in this article works, with slight complications; the limit is less surprising. We may want to let the value of $\mu_0$ ($=$ the constant value of $\mu$ within the square) vary on all of $\D$ and even beyond: it seems to be possible to do this; one should find an interpretation for it. We may let $\mu$ tend to $1$ along directions not contained in $\R$ and study how the limit depends on the direction.

The fact that when a finite type affine surface can be uniformized onto the Riemann sphere $\S$ (finitely punctured) the affine coordinates on $\S$ are given by a generalized Schwarz-Christoffel formula was probably already known, but is not so easy to find in the literature.\footnote{It can be found in a slightly different form, like in \cite{Gil}, according to whom it was already known from Schwarz and Schläfli.}
Concerning the idea of using affine surfaces to uniformize Beltrami forms that are piecewise constant on polygons, the author is not aware of this having already been tried before. This should be studied further.

Concerning affine surfaces and the distortion derivative: if a holomorphic function $\eta$ has a singularity of polar type of order $>1$, can we find and classify models of the affine surface structure that $\eta$ defines near the singularity? Can we do the same for projective structures, i.e.\ with the Schwarzian derivative replacing the distortion derivative?

\end{document}